\documentclass[12pt,a4paper,twoside]{article}

\usepackage{stmaryrd}

\usepackage{hyperref}

\newcommand{\pageformat}[6]{\setlength{\hoffset}{-1in}
                  \setlength{\voffset}{-1in}
                  \addtolength{\hoffset}{#5}
                            \addtolength{\voffset}{#6}
                            \setlength{\oddsidemargin}{#1}
                            \setlength{\evensidemargin}{#2}
                            \setlength{\textwidth}{\paperwidth}
                  \addtolength{\textwidth}{-\oddsidemargin}
                  \addtolength{\textwidth}{-\evensidemargin}
                  \addtolength{\textwidth}{-\marginparsep}
                  \addtolength{\textwidth}{-\marginparwidth}
                            \setlength{\topmargin}{#3}
                            \setlength{\textheight}{\paperheight}
                  \addtolength{\textheight}{-\topmargin}
                  \addtolength{\textheight}{-\headheight}
                  \addtolength{\textheight}{-\headsep}
                  \addtolength{\textheight}{-\footskip}
                  \addtolength{\textheight}{-#4}}
\pageformat{2cm}{3cm}{25mm}{25mm}{1pt}{0pt}

\usepackage{ifthen}
\newboolean{article}
    \setboolean{article}{true}
\newboolean{report}
\newboolean{book}
\newboolean{letter}
\newboolean{german}
\newboolean{italian}
\newboolean{nobaselinestretch}
\newboolean{nosectionappendix}
\newboolean{oldtoc}
\newboolean{nosectionequn}
\newboolean{notheorem}

\ifthenelse{\boolean{german}}{
    \usepackage{german}}{}

\usepackage[latin1]{inputenc}

\makeatletter
\ifthenelse{\boolean{article}\and\not\boolean{oldtoc}}{
    \renewcommand*\l@section[2]{%
        \ifnum \c@tocdepth >\z@
            \addpenalty\@secpenalty
            \addvspace{0.0em \@plus\p@}%
            \setlength\@tempdima{1.5em}%
            \begingroup
                \parindent \z@ \rightskip \@pnumwidth
                \parfillskip -\@pnumwidth
                \leavevmode 
                \advance\leftskip\@tempdima
                \hskip -\leftskip
                #1\nobreak\dotfill \nobreak\hb@xt@\@pnumwidth{\hss #2}\par
            \endgroup
        \fi}}{}
\makeatother

\usepackage{amsmath}
\usepackage{amssymb}
\usepackage[mathscr]{eucal}

\ifthenelse{\boolean{notheorem}}{}{
    \usepackage{theorem}}

\ifthenelse{\boolean{nobaselinestretch}}{}{
    \renewcommand{\baselinestretch}{1.25}}

\newenvironment{env}[2]{\begin{#1}#2\end{#1}}{}
    \newcommand{\beq}[1]{\begin{env}{equation}{#1}}
    \newcommand{\beqn}[1]{\begin{env}{equation*}{#1}}
    \newcommand{\bal}[1]{\begin{env}{align}{#1}}
    \newcommand{\baln}[1]{\begin{env}{align*}{#1}}
    \newcommand{\bga}[1]{\begin{env}{gather}{#1}}
    \newcommand{\bgan}[1]{\begin{env}{gather*}{#1}}
    \newcommand{\bflal}[1]{\begin{env}{flalign}{#1}}
    \newcommand{\bflaln}[1]{\begin{env}{flalign*}{#1}}
    \newcommand{\bmu}[1]{\begin{env}{multline}{#1}}
    \newcommand{\bmun}[1]{\begin{env}{multline*}{#1}}
    \newcommand{\bsp}[1]{\begin{env}{split}{#1}}

    \newcommand{\eeq}{\end{env}}
    \newcommand{\eeqn}{\end{env}}
    \newcommand{\eal}{\end{env}}
    \newcommand{\ealn}{\end{env}}
    \newcommand{\ega}{\end{env}}
    \newcommand{\egan}{\end{env}}
    \newcommand{\eflal}{\end{env}}
    \newcommand{\eflaln}{\end{env}}
    \newcommand{\emu}{\end{env}}
    \newcommand{\emun}{\end{env}}
    \newcommand{\esp}{\end{env}}

\newcommand{\lf}{\vspace{2ex}}
\newcommand{\bulletline}[1][]{\lf\noindent~\hfill$\bullet\bullet\bullet$\hfill~

\lf\noindent\bf{#1}}

\renewcommand{\bf}[1]{\textbf{#1}}
\renewcommand{\it}[1]{\textit{#1}}

\renewcommand{\sf}[1]{\textsf{#1}}

\renewcommand{\tt}[1]{\texttt{#1}}
\newcommand{\hl}[1]{\bf{\it{#1}}}

\newcommand{\mbf}[1]{\mathbf{#1}}
\newcommand{\msf}[1]{\text{\small $\sf{#1}$}}

\newcommand{\cmc}[1]{\mathcal{#1}}
\newcommand{\eus}[1]{\mathscr{#1}}
\newcommand{\euf}[1]{\mathfrak{#1}}
\newcommand{\bb}[1]{\mathbb{#1}}

\newcommand{\mtiny}[1]{{\setlength{\arraycolsep}{.3ex}\text{\tiny$#1$}}}
\newcommand{\nbd}[1]{$#1$\nobreakdash--}
\newcommand{\ol}[1]{\overline{#1}}

\newcommand{\wt}[1]{\widetilde{#1}}
\newcommand{\wh}[1]{\widehat{#1}}

\newcommand{\vt}{\vartheta}

\newcommand{\norm}[1]{\left\lVert#1\right\rVert}

\newcommand{\snorm}[1]{\norm{\smash{#1}}}

\newcommand{\bfam}[1]{\bigl(#1\bigr)}
\newcommand{\Bfam}[1]{\Bigl(#1\Bigr)}
\newcommand{\AB}[1]{\langle#1\rangle}

\newcommand{\CB}[1]{\{#1\}}
\newcommand{\bCB}[1]{\bigl\{#1\bigr\}}
\newcommand{\BCB}[1]{\Bigl\{#1\Bigr\}}
\newcommand{\SB}[1]{[#1]}
\newcommand{\bSB}[1]{\bigl[#1\bigr]}

\newcommand{\Matrix}[1]{\begin{pmatrix}#1\end{pmatrix}}

\newcommand{\tMatrix}[1]{\mtiny{\Matrix{#1}}}

\newcommand{\rtMatrix}[1]{\raisebox{.3ex}{\tMatrix{#1}}}

\newcommand{\sbars}[1]{\:\bar{#1}^s\:}
\newcommand{\sodots}{\sbars{\odot}}
\newcommand{\set}[2][]{
    \ifthenelse{\equal{#1}{}}{
        \CB{#2}}{
        \CB{#1~|~#2}}}
\newcommand{\bset}[2][]{
    \ifthenelse{\equal{#1}{}}{
        \bCB{#2}}{
        \bCB{#1~|~#2}}}
\newcommand{\Bset}[2][]{
    \ifthenelse{\equal{#1}{}}{
        \BCB{#2}}{
        \BCB{#1~\big|~#2}}}
\newcommand{\zero}{\CB{0}}

\DeclareMathOperator{\ls}{\normalfont\msf{span}}

\DeclareMathOperator{\cls}{\ol{\ls}}

\DeclareMathOperator{\id}{\normalfont\msf{id}}

\newcommand{\C}{\bb{C}}

\newcommand{\E}{\bb{E}}

\newcommand{\N}{\bb{N}}

\newcommand{\R}{\bb{R}}

\newcommand{\T}{\bb{T}}

\newcommand{\cA}{\cmc{A}}
\newcommand{\cB}{\cmc{B}}
\newcommand{\cC}{\cmc{C}}

\newcommand{\sB}{\eus{B}}

\newcommand{\sF}{\eus{F}}

\newcommand{\sS}{\eus{S}}

\newcommand{\en}{\euf{n}}

\newcommand{\eu}{\euf{u}}

\newcommand{\eB}{\euf{B}}

\newcommand{\eH}{\euf{H}}

\newcommand{\U}{\mbf{1}}

\newcommand{\I}{{I\!\!\!\;I}}

\ifthenelse{\boolean{nosectionequn}}{}{
    \numberwithin{equation}{section}
    }

\ifthenelse{\boolean{article}\or\boolean{letter}\or\boolean{nosectionequn}}{
    \setboolean{nosectionappendix}{true}}{}
\ifthenelse{\boolean{nosectionappendix}}{}{
    \renewcommand{\appendix}{
        \chapter*{\appendixname}
        \addcontentsline{toc}{chapter}{\appendixname}
        \renewcommand{\thesection}{\Alph{section}}
        \setcounter{section}{0}}}
   
\ifthenelse{\boolean{report}\or\boolean{book}}{
    }{}

\ifthenelse{\boolean{notheorem}}{}{
        \newcommand{\mnname}{Mathematical note.}
        \newcommand{\enname}{End of the note.}
        \newcommand{\definame}{Definition.}
        \newcommand{\propname}{Proposition.}
        \newcommand{\lemname}{Lemma.}
        \newcommand{\exname}{Example.}
        \newcommand{\exername}{Exercise.}
        \newcommand{\remname}{Remark.}
        \newcommand{\obname}{Observation.}
        \newcommand{\thmname}{Theorem.}
        \newcommand{\corname}{Corollary.}
        \newcommand{\proofname}{Proof.}
    \ifthenelse{\boolean{german}}{
        \renewcommand{\mnname}{Mathematische Notiz.}
        \renewcommand{\enname}{Ende der Notiz.}
        \renewcommand{\exname}{Beispiel.}
        \renewcommand{\exername}{ï¿½bung.}
        \renewcommand{\remname}{Bemerkung.}
        \renewcommand{\obname}{Beobachtung.}
        \renewcommand{\thmname}{Satz.}
        \renewcommand{\corname}{Korollar.}
        \renewcommand{\proofname}{Beweis.}}{}
    \ifthenelse{\boolean{italian}}{
        \renewcommand{\mnname}{Nota matematica.}
        \renewcommand{\enname}{Fina della nota.}
        \renewcommand{\definame}{Definizione.}
        \renewcommand{\propname}{Proposizione.}
        \renewcommand{\exname}{Esempio.}
        \renewcommand{\exername}{Esercizio.}
        \renewcommand{\remname}{Nota.}
        \renewcommand{\obname}{Osservazione.}
        \renewcommand{\thmname}{Teorema.}
        \renewcommand{\corname}{Corollario.}
        \renewcommand{\proofname}{Dimostrazione.}

       \renewcommand{\appendixname}{Appendice}

       }{}
    \theoremheaderfont{\normalfont\bfseries}
    \theoremstyle{change}
        \theorembodyfont{\rmfamily}
            \newtheorem{emp}{}[section]
                \newcommand{\bemp}[1][]{
                    \begin{emp}\hskip-\labelsep\bf{#1}\hskip\labelsep}
                \newcommand{\eemp}{\end{emp}}
\newtheorem{itemp}[emp]{}
                \newcommand{\bitemp}[1][]{
                    \begin{itemp}\hskip-\labelsep\bf{#1}\hskip\labelsep\normalfont\itshape}
                \newcommand{\eitemp}{\end{itemp}}
            \newtheorem{mn}[emp]{\mnname}
                \newcommand{\bnm}{\begin{mn}~\begin{quotation}\renewcommand{\baselinestretch}{1}\small\noindent\ignorespaces}
                \newcommand{\enm}{\end{quotation}\hfill\bf{\enname}\end{mn}}
            \newtheorem{ex}[emp]{\exname}
                \newcommand{\bex}{\begin{ex}}
                \newcommand{\eex}{\end{ex}}
            \newtheorem{exer}[emp]{\exername}
                \newcommand{\bexer}{\begin{exer}}
                \newcommand{\eexer}{\end{exer}}
            \newtheorem{defi}[emp]{\definame}
                \newcommand{\bdefi}{\begin{defi}}
                \newcommand{\edefi}{\end{defi}}
            \newtheorem{rem}[emp]{\remname}
                \newcommand{\brem}{\begin{rem}}
                \newcommand{\erem}{\end{rem}}
            \newtheorem{ob}[emp]{\obname}
                \newcommand{\bob}{\begin{ob}}
                \newcommand{\eob}{\end{ob}}
        \theorembodyfont{\normalfont\itshape}
            \newtheorem{thm}[emp]{\thmname}
                \newcommand{\bthm}{\begin{thm}}
                \newcommand{\ethm}{\end{thm}}
            \newtheorem{prop}[emp]{\propname}
                \newcommand{\bprop}{\begin{prop}}
                \newcommand{\eprop}{\end{prop}}
            \newtheorem{cor}[emp]{\corname}
                \newcommand{\bcor}{\begin{cor}}
                \newcommand{\ecor}{\end{cor}}
            \newtheorem{lem}[emp]{\lemname}
                \newcommand{\blem}{\begin{lem}}
                \newcommand{\elem}{\end{lem}}
\newenvironment{empn}[1]{\lf\noindent\bf{#1}\ignorespaces\hskip\labelsep}{\lf}
		\newcommand{\bempn}[1]{\begin{empn}{#1}}
		\newcommand{\eempn}{\end{empn}}
		\newcommand{\bitempn}[1]{\begin{empn}{#1}\normalfont\itshape}
		\newcommand{\eitempn}{\end{empn}}
                \newcommand{\bnmn}{\begin{empn}{\mnname}~\begin{quotation}\renewcommand{\baselinestretch}{1}\small\noindent\ignorespaces}
                \newcommand{\enmn}{\end{quotation}\hfill\bf{\enname}\end{empn}}
		\newcommand{\bexn}{\begin{empn}{\exname}}
		\newcommand{\eexn}{\end{empn}}
		\newcommand{\bexern}{\begin{empn}{\exername}}
		\newcommand{\eexern}{\end{empn}}
		\newcommand{\bdefin}{\begin{empn}{\definame}}
		\newcommand{\edefin}{\end{empn}}
		\newcommand{\bremn}{\begin{empn}{\remname}}
		\newcommand{\eremn}{\end{empn}}
		\newcommand{\bobn}{\begin{empn}{\obname}}
		\newcommand{\eobn}{\end{empn}}

\newcommand{\qedsymbol}{~\rule[-0.35mm]{2mm}{2mm}}
    \newcounter{proof}[emp]
    \newenvironment{Proof}[1]{
        \vspace{1ex}
        \renewcommand{\item}[1][\stepcounter{proof}(\roman{proof})]%
            {##1\hskip\labelsep}
        \noindent\textsc{#1\hskip\labelsep}}{
        \nolinebreak\qedsymbol}
    \newcommand{\proof}[1][\proofname]{
        \begin{Proof}{#1}\ignorespaces}
    \newcommand{\qed}{\end{Proof}}
    \newcommand{\noqed}{
        \renewcommand{\qedsymbol}{}
        \end{Proof}}}
    \ifthenelse{\boolean{italian}}{
        \renewcommand{\proofname}{Dimostrazione.}}{}

\usepackage[varg]{txfonts}
\usepackage{bbm}

\usepackage[matrix,arrow,curve]{xy}

\renewcommand{\thefootnote}{[\alph{footnote}]}

\begin{document}



\renewcommand{\thefootnote}{{(\alph{footnote})}}


\bibliographystyle{amsalpha}
\newcommand{\ins}[1]{}

\title{Paired $E_0$--Semigroups}
\author{}
\author{Michael Skeide}
\date{}

{
\renewcommand{\baselinestretch}{1}
\maketitle


\begin{abstract}
\noindent
In these notes we prove two main results:

1) It is well-known \ins{(Arveson \cite{Arv89}) }that two strongly continuous \nbd{E_0}semigroups on $\sB(H)$ can be \it{paired} \ins{in the sense of Powers and Robinson \cite{PoRo89} }if and only if they have anti-isomorphic Arveson systems. For a new notion of \it{pairing} (which coincides only for $\sB(H)$ with the existing one), we show: For a von Neumann algebra $\cB$, a strongly continuous \nbd{E_0}semigroup $\vt$ on $\cB$ and a strongly continuous \nbd{E_0}semigroup $\vt'$ on $\cB'$ can be \it{paired}
if and only if their product systems are commutants of each other.

2) On the way to prove the former, \it{en passant} we have to fill in a long standing important gap in the theory of intertwiner product systems \it{à la} Arveson: \it{Intertwiner product systems} of faithful strongly continuous \nbd{E_0}semigroups on von Neumann algebras have \it{enough strongly continuous sections}.

We explain why both results are entirely out of reach for Arveson's methods \cite{Arv89,Arv90} and depend essentially on the alternative approach from Skeide \cite{Ske16}.
\end{abstract}



\tableofcontents
}





\section{Introduction}

Throughout, $\cB\subset\sB(G)$ is a von Neumann algebra acting (nondegenerately!)\ on a Hilbert space $G$, and $\cB'$ its commutant. An \hl{\nbd{E_0}semigroup} on $\cB$ is a semigroup $\bfam{\vt_t}_{t\in\R_+}$ of normal unital endomorphisms of $\cB$. Endomorphism means \nbd{*}endomorphism. Apart from normality of each $\vt_t$, no continuity with $t$ is required. (For more conventions and notations, see the end of this introduction.)

\lf
\nbd{E_0}semigroups arise naturally by restricting automorphism groups on $\sB(G)$ to \it{future} and \it{past} algebras invariant for negative times and positive times, respectively. It is known from Arveson and Kishimoto \cite{ArKi92} that any faithful \nbd{E_0}semigroup on a \nbd{W^*}algebra may be extended to an automorphism group embedding the \nbd{W^*}algebra suitably as a von Neumann algebra in some $\sB(G)$. (See also the proof of \cite[Theorem B.36]{Ske16} or, for the \nbd{C^*}case, \cite{Ske11a}.%
\footnote{
For $\cB=\sB(G)$, the result occurred first as \cite[Corollary 5.21]{Arv90}. Like our proofs for general $\cB$ in \cite{Ske11a,Ske16}, the proof in \cite{Arv90} depends on the result that every faithful product system admits, in our terminology, a right dilation (see Definition \ref{rdthm}), that is, a faithful representation (see Definition \ref{PSrepdefi}). Once this \it{fundamental theorem on Arveson systems} is established, the construction is easy. The proof of the fundamental theorem on Arveson systems is the main result of \cite{Arv90}; its proof there is deep and hard. But there are simpler proofs in \cite{Ske06,Arv06}, and a combination of both these proofs \cite{Ske06a} generalizes to more general (faithful) product systems than Arveson systems.

The proof in \cite{ArKi92} is totally different and of intermediate difficulty.
}%
) By Proposition \ref{bashisprop} below, this situation is closely related to our new notion of \it{pairing} \nbd{E_0}semigroups and, for \nbd{E_0}semigroups on $\sB(H)$, to \it{pairing} as considered by Powers and Robinson \cite{PoRo89} (see Remark \ref{PRArem} below). Powers \cite{Pow88,Pow87} started studying \nbd{E_0}semi\-groups on $\sB(H)$ up to \it{cocycle conjugacy}. Arveson \cite{Arv89} classified them up to cocycle conjugacy by associating with each \nbd{E_0}semigroup on $\sB(H)$ a \it{product system of Hilbert spaces} (a so-called \it{Arveson system}) and showing that the isomorphism classes of Arveson systems are in one-to-one correspondence with the cocycle conjugacy classes of \nbd{E_0}semigroups on $\sB(H)$. The construction of product systems from quantum dynamics such as \nbd{E_0}semigroups has been modified and/or generalized in various ways; \cite{Bha96,BhSk00,MuSo02,Fow02,Ske03c,Ale04,Ske06,Arv06,MSS06,Ske09,Ske16,MaSr17,ShaSk23}.%
\footnote{
Arveson's \it{intertwiner construction} \cite{Arv89} generalizes to what we call here the \it{intertwiner system} (\it{à la} Arveson) of an \nbd{E_0}semigroup on a general von Neumann algebra as done in Skeide \cite{Ske03c}. It coincides with Alevras' construction \cite{Ale04} for type II$_1$ factors, but is not limited to these. The von Neumann correspondences arising in the latter, probably occurred already in Alevras' PhD-thesis 1995, but despite several attempts we never managed to get our hands on a copy, so we could not verify. The first occurrence of product systems of correspondences in an article seems to be in Bhat and Skeide \cite{BhSk00} in the context of dilations of Markov semigroups. Bhat \cite{Bha96} (see Theorem \ref{Bthm} below) provided a different construction of an Arveson system for \nbd{E_0}semigroups on $\sB(H)$ (which is effectively anti-isomorphic to Arveson's), which has been generalized in Skeide \cite{Ske02} (see Theorem \ref{Sthm} below) to \nbd{E_0}semigroups on $\sB^a(E)$. Fowler \cite{Fow02} was the first one to consider product systems indexed by general monoids, while Shalit and Skeide \cite{ShaSk23} push this forward to the extreme (also in relation with structures more general than product systems).

We may very well add to the literature list \cite{Ara70,PaSchm72}, \cite{MSchue93}, and \cite{HKK04p} for the first occurrences of product, sub-, and superproduct systems, respectively, (without recognizing the structures as such), and \cite{ShaSo09,BhMu10} and \cite{MaSr13} for the first explicit definitions of sub- and superproduct systems, respectively.
}
(See also Remark \ref{rem3rd} below, for a possible extension of the present work to open systems with invariant faithful states.)

Arveson himself showed (see \cite[Theorem 3.5.5]{Arv03}) that two strongly continuous $E_0$-- semi\-groups on $\sB(H)$ can be paired (in the sense of \cite{PoRo89}; see Remark \ref{PRArem}) if and only if they have anti-isomorphic Arveson systems. Here, we propose a new definition of pairing (coinciding only for $\sB(H)$ with the definition in \cite{PoRo89} for general von Neumann algebras), and show a result that generalizes Arveson's provided we take into account the fact (mentioned in  the places in the surveys Skeide \cite{Ske05a,Ske08a} indicated in Remark \ref{opprem} below) that the commutant of an Arveson system is an Arveson system anti-isomorphic to the original one.

\bdefi \label{paireddefi}
A \hl{pairing} for an \nbd{E_0}semigroup $\vt$ on $\cB$ and an \nbd{E_0}semigroup $\vt'$ on $\cB'$ is an automorphism group $\alpha$ on $\sB(G)$ such that $\vt_t(b)=\alpha_{-t}(b)$ $(t\ge0)$ and $\vt'_t(b')=\alpha_t(b')$ $(t\ge0)$ for all $b\in\cB$ and $b'\in\cB'$. We also say, $\vt$ and $\vt'$ are \hl{paired} (via $\alpha$).
\edefi

In Section \ref{cpairSEC} we will prove the following result.

\bthm \label{cpairedthm}
A strongly continuous \nbd{E_0}semigroup $\vt$ on $\cB$ and a strongly continuous \nbd{E_0}semi\-group $\vt'$ on $\cB'$ can be paired via a strongly continuous $\alpha$ if and only if their strongly continuous product systems are strongly continuous commutants of each other.
\ethm

We will prove this as a consequence of Theorem \ref{pairedthm} below without continuity conditions, plus various results from Skeide \cite[Appendix B]{Ske16} incorporating continuity, plus the known form of strongly continuous automorphisms groups on $\sB(G)$ going back to Wigner \cite{Wig39}, plus the new result Theorem \ref{ecthm} that  answers a long standing question about so-called \it{intertwiner product systems à la Arveson} of strongly continuous faithful \nbd{E_0}semigroups in the affirmative:

\bthm \label{ecthm}
Intertwiner product systems of faithful strongly continuous \nbd{E_0}semigroups have enough strongly continuous sections.
\ethm

While the main motivation for these notes is proving Theorem \ref{cpairedthm}, for the general theory of product systems the other main result, Theorem \ref{ecthm}, might very well be considered more important. Note that Theorem \ref{ecthm} is new even for \nbd{E_0}semigroups on $\sB(H)$ for separable Hilbert spaces $H$. (In this case, Arveson \cite[Lemma 2.3]{Arv89} has a partial version of Theorem \ref{ecthm}, namely, that the intertwiner spaces admit (enough) sections that act strongly continuously on $H$ for $t>0$;%
\footnote{ \label{AnoncFN}
\cite[Lemma 2.3]{Arv89} states that for a strongly continuous non-automorphic \nbd{E_0}semigroup $\vt$ on $\sB(H)$ ($H$ infinite-dimensional and separable) and any fixed $t_0>0$ there is a strongly continuous family (\bf{not} semigroup!) $\bfam{U_t}_{t>0}$ of unitaries $U_t\in\sB(H)$ such that $\vt_t(a)=U_t\vt_{t_0}(a)U_t^*$ for all $t>0$. Clearly, this means that if $a$ intertwines $\vt_{t_0}$ (that is, if $\vt_{t_0}(b)a=ab$ for all $b\in\sB(H)$), then $U_ta$ intertwines $\vt_t$. So for every $c$ intertwining $\vt_t$ for some $t>0$, the section $\bfam{U_sU_t^*c}_{s>0}$ is a strongly continuous section of intertwiners for $s>0$ and hits $c$ for $s=t$. But there is no indication that the section could be continued to $s=0$. On the contrary, the proof of \cite[Lemma 2.3]{Arv89} is based on Dixmier \cite[Lemma 10.8.7]{Dix77} on the structure of continuous fields of Hilbert spaces with constant dimension $\aleph_0$, and the intertwiner space of $\vt$ at $t=0$ is one-dimensional.

Clearly, the unitary family $\bfam{U_t}_{t>0}$ cannot be continued to $t=0$ by a unitary $U_0$, giving back $\vt_0=\id_{\sB(H)}$ as $U_0\vt_{t_0}(\bullet)U_0^*$ because, then, all $\vt_t$ would be automorphisms. So, there is no way to get strongly continuous sections of intertwiners by means of something like \cite[Lemma 2.3]{Arv89} continued to $t=0$.
}
but the information -- crucial for proving Theorem \ref{cpairedthm} -- that there are such sections continuous also at $t=0$ (necessarily taking there a value in the one-dimensional space $\id_H\C$) is entirely missing. Other authors who compute intertwiner product systems (such as  Alevras \cite{Ale04} or Margetts and Srinivasan \cite{MaSr13}) from \nbd{E_0}semigroups on type II$_1$ factors address only Borel isomorphisms between the trivial bundle (obviously strongly continuous) and the product system for $t>0$, but do not address the question if the Borel isomorphism can be chosen strongly continuous at least for $t>0$.)

Another prerequisite for the proof of Theorem \ref{cpairedthm} is the following algebraic version:

\bthm \label{pairedthm}
An \nbd{E_0}semigroup $\vt$ on $\cB$ and an \nbd{E_0}semigroup $\vt'$ on $\cB'$ can be paired if and only if their product systems are commutants of each other modulo a multiplier over $\R_+$.
\ethm

The missing terminology necessary to give a precise meaning to Theorems \ref{pairedthm} and \ref{ecthm} (especially, the meaning of Theorem \ref{ecthm} is clear to everybody who ever has constructed an Arveson system \it{à la} Arveson: the intertwiner product system has  \hl{enough strongly continuous sections} if the strongly continuous sections of intertwiners hit each point in the intertwiner product system; see also Definition \ref{scPSdthm}) is easy to explain; this will be done still in this introduction. On the contrary, for just interpreting the strongly continuous version in Theorem \ref{cpairedthm}, we have to explain terminology and results about strongly continuous product systems and their strongly continuous commutants from Skeide \cite[Section 12+Appendix B]{Ske16}; the latter will be done not before Section \ref{E0psSEC}. While the meaning of Theorem \ref{ecthm} is, as we just mentioned, immediate, its proof in Section \ref{E0psSEC} is not and requires fully the notions and results from that section. (For that proof, in Section \ref{E0psSEC} we have to fix certain isomorphisms with a precision never required before; the reason, why in earlier papers we could be sloppier and work ``up to isomorphism'', is explained in Footnote \ref{rduniFN}.)

As promised, we now explain immediately the missing terminology for the algebraic Theorem \ref{pairedthm}. (For the reader who does not yet know the definition of product system, this can be looked up in Theorem and Definition \ref{Sthm}.) And after having introduced this terminology and some add-ons about multipliers in Section \ref{1dimAsSEC}, the reader interested only in Theorem \ref{pairedthm} may pass directly to its proof in Section \ref{pairSEC}. It seems, however, appropriate to mention already now that the product systems in the second item of the following are the same used also in the proof of Theorem \ref{cpairedthm}; just their strongly continuous structure is added later. Also the identifications used in the proof of Theorem \ref{pairedthm} in Section \ref{pairSEC} are the same as those for the proof of Theorem \ref{cpairedthm} in Section \ref{cpairSEC}.

\begin{itemize}
\item
A \hl{multiplier} over $\R_+$ is a function $m\colon\R_+^2\rightarrow\T$ satisfying Equation \eqref{mult}.  Multipliers determine the structure of one-dimensional \it{Arveson systems} (that is, product systems of one-dimensional Hilbert spaces) -- and with that, they determine the structure of not necessarily continuous automorphism (semi)groups on $\sB(G)$. (We explain this in Section \ref{1dimAsSEC}.) Assuming the reader knows algebraic \it{product systems} (that is, no technical conditions such as continuity or measurability; see again Theorem and Definition \ref{Sthm}) and related notions, a multiplier $m$ may act on a product system $E^\odot=\bfam{E_t}_{t\in\R_+}$ by replacing the product $(x_s,y_t)\mapsto u_{s,t}(x_s\odot y_t)$ with the product $(x_s,y_t)\mapsto u_{s,t}(x_s\odot y_t)m(s,t)$, equipping the same family with another product system structure. $E^\odot$ is isomorphic to $F^\odot$ \hl{via the multiplier $m$} if $E^\odot$ with the new product induced by $m$ is isomorphic to $F^\odot$.

\item
There are various methods to associate with an \nbd{E_0}semigroup a product system of von Neumann correspondences. (It is noteworthy that Theorem \ref{pairedthm} does not depend on the choice of the method, provided we choose for $\vt$ and $\vt'$ the same method.) Most methods associate with an \nbd{E_0}semigroup $\vt$ on a von Neumann algebra $\cB\subset\sB(G)$ (or, more generally, on $\sB^a(E)$, where $E$ is a von Neumann \nbd{\cB}module) a product system of von Neumann \nbd{\cB}correspondences, while the intertwiner product system \it{à la} Arveson consists of von Neumann \nbd{\cB'}correspondences. The former are all isomorphic among themselves, while the latter is isomorphic to the commutant of the former.

For our purposes here (that is, looking only at \nbd{E_0}semigroups $\vt$ on $\cB=\sB^a(\cB)\subset\sB(G)$), what we are going to call \bf{the} product system of $\vt$ is particularly simple: \bf{The} \hl{product system $E^\odot=\bfam{E_t}_{t\in\R_+}$ of $\vt$} is given as $E_t:={_t}\cB$ (which is shorthand for ${_{\vt_t}}\cB$), that is, $\AB{x_t,y_t}:=x_t^*y_t$ (with the only possible right multiplication making it an inner product, namely, operator multiplication in $\cB\subset\sB(G)$) and with left multiplication $b.x_t:=\vt_t(b)x_t$. And the product of the product system is $u_{s,t}(x_s\odot y_t):=\vt_t(x_s)y_t$. (This is compatible with the more general situation, when $\vt$ acts on $\sB^a(E)$ for a von Neumann \nbd{\cB}module $E$ with a unit vector $\xi$ as the special case $E=\cB$ and $\xi=\U$; see Observation \ref{TPSob}.) The \hl{intertwiner product system} \it{à la} Arveson (or simply \hl{intertwiner system}) $F'^\odot=\bfam{F'_t}_{t\in\R_+}$ of $\vt$ is given by the von Neumann \nbd{\cB'}correspondences%
\footnote{ \label{CBMdefFN}
For any \nbd{\cB}bimodule $M$ we define its \hl{\nbd{\cB}center} $C_\cB(M):=\CB{x\in M\colon bx=xb~(b\in\cB)}$.
}
 \beqn{
F'_t
~:=~
C_\cB(\sB(G,{_t}G))
~=~
\BCB{x'_t\in\sB(G)\colon\vt_t(b)x'_t=x'_tb~(b\in\cB)}
}\eeqn
with \nbd{\cB'}bimodule structure by operator multiplication and inner product $\AB{x'_t,y'_t}:={x'_t}^*y'_t$ and with product system structure $u'_{t,s}(x'_t\odot y'_s):=x'_ty'_s$ (also operator multiplication).
\end{itemize}

\noindent
In our simplified setting ($\sB^a(\cB)=\cB$ instead of $\sB^a(E)$ for general $E$), $F'^\odot$ is not only (canonically) isomorphic to the commutant of $E^\odot$; it really \bf{is} the commutant of $E^\odot$. So, taking in this statement (it can be seen easily by limiting the discussion in Skeide \cite[Section 2]{Ske03c} to our case; but we also prove it in Observation \ref{commob}), we are now already ready to interpret Theorem \ref{pairedthm} and, after adding the necessities about multipliers in Section \ref{1dimAsSEC}, we are also able to prove Theorem \ref{pairedthm}: Take \bf{the} product system of $\vt$, $E^\odot$, and the commutant $F^\odot:=(F')'^\odot$ of \bf{the} product system of $\vt'$, $F'^\odot$, (namely, the intertwiner system of $\vt'$) and show that it is possible to construct an isomorphism modulo a multiplier $m$ between $E^\odot$ and $F^\odot$ if and only if $\vt$ and $\vt'$ can be paired.

We repeat: The product systems that occur in Theorem \ref{cpairedthm} and the isomorphisms constructed among them are the same as for Theorem \ref{pairedthm}. It is one benefit in the strongly continuous case that, thanks to Wigner's theorem \cite{Wig39}, the multiplier goes away (leading to the simplified version of Theorem \ref{pairedthm} where $\alpha$ is inner, Theorem \ref{ipairedthm}). However, all the rest is considerably more complicated. $E^\odot$ and $F^\odot$ have to be equipped with a strongly continuous structure (as defined in Skeide \cite[Section 12]{Ske16}). For $E^\odot$ this is easy and can be done directly: The strongly continuous sections are just those functions $t\mapsto x_t\in E_t={_t}\cB=\cB\subset\sB(G)$ that are strongly continuous maps into $\sB(G)$. But the same does not work directly for the intertwiner product system of $\vt'$, $F^\odot$. We simply do not \it{a priori} know if there are enough strongly continuous sections of intertwiners; the result that tells us that this is true, Theorem \ref{ecthm}, has not been dealt with in \cite{Ske16}. And in fact, only after recognizing $F^\odot$ as the strongly continuous commutant (in the sense of \cite[Appendix B]{Ske16}) of \bf{the} strongly continuous product system of $\vt'$, $F'^\odot$, we get the right candidates for the strongly continuous sections of $F^\odot$ and can show that they, indeed, act, when considered sections of the intertwiner system of $\vt'$, strongly continuously on $G$.

Note that for just defining the strongly continuous commutant of a strongly continuous product system, even of a product system as simple as our $E^\odot$ here, we need the full construction of \bf{the} strongly continuous product system of an \nbd{E_0}semigroup acting on $\sB^a(E)$ based on a unit vector $\xi\in E$; see Point \ref{outob} in the beginning of Section  \ref{E0psSEC}. Also the identification of the intertwiner system with the commutant system is required explicitly -- or at least more so, than we ever did in our earlier papers \cite{Ske03c,Ske08,Ske09,Ske16} (see also Footnote \ref{rduniFN}). For this, it is unavoidable to work in the framework of \bf{concrete} von Neumann correspondences (Skeide \cite{Ske06b}) in which the commutant is defined as a bijective functor and not only as a natural equivalence. We give the necessary details in Section \ref{E0psSEC}.

\bulletline
The following three remarks -- apart from giving additional detail, also explaining why our approach here is very different from Arveson's, which does, for principal reasons, not generalize to arbitrary von Neumann algebras --, in particular, point into the direction of future developments based on the present notes and/or generalizing them. It is safe to skip these remarks.

~

~

\brem \label{rem1st}
Let us have a quick word on ``\it{faithful}''. Clearly, an \nbd{E_0}semigroup $\vt$, in order to possess an automorphic extension, has to be faithful.%
\footnote{
A warning: This hypothesis is missing in \cite[Theorem A]{ArKi92}
}
Faithfulness of $\vt$ corresponds to faithfulness of \bf{the} product system of $\vt$. And, under commutant, faithfulness of a product system corresponds to strong fullness of its commutant, and \it{vice versa}. Recall that \bf{the} product system of an \nbd{E_0}semigroup is necessarily strongly full. So, speaking about pairings, we are necessarily speaking about product systems that are strongly full and faithful; and this class is invariant under commutant. Moreover, in the strongly continuous case, the construction of the strongly continuous commutant of a strongly continuous product system in Skeide \cite[Appendix B]{Ske16} is limited to strongly full and faithful product systems. This is the reason, why we get Theorem \ref{ecthm} only for faithful \nbd{E_0}semigroups. (This, clearly, includes all \nbd{E_0}semigroups on factors.)

Missing unitality of $\vt$ (a so-called \nbd{E}semigroup) may be easily discussed away by passing to the unitalization $\wt{\vt}$ (acting on the unitalization $\wt{\cB}=\C\oplus\cB\subset\sB(\C\oplus G)$ of $\cB$%
\footnote{
Observe, however, that this is not possible in a theory limited to $\cB=\sB(G)$, because $\wt{\sB(G)}$ is not another $\sB(\wt{G})$, but a proper subalgebra not isomorphic to any type I factor; in fact, it has nontrivial center.
}%
). We do not know how to do a ``\it{faithfulization}'' for an \nbd{E_0}semigroup.

We do believe that the theory of strongly full and faithful strongly continuous product systems and their commutants from \cite{Ske16} (based on the construction of an \nbd{E_0}semigroup for any strongly full strongly continuous product system) can be pushed forward to a theory of arbitrary strongly continuous product systems and their commutants (to be based on the well-known construction of an \nbd{E}semigroup by letting act $E^\odot$ from the right on the direct integral $\int^\oplus E_t\,dt$). The algebraic idea behind the construction of this \nbd{E}semigroup is much simpler than that (in the strongly full case) of an \nbd{E_0}semigroup, but the technical problems will be different. A complete discussion would include a variant of Theorem \ref{ecthm} for the intertwiner systems of arbitrary \nbd{E}semigroups; but this would be another paper.
\erem

\brem \label{rem2nd}
Theorem \ref{ecthm} says that the intertwiner system of a (faithful) strongly continuous \nbd{E_0}semi\-group on $\cB\subset\sB(G)$ is a strongly continuous subbundle of the trivial strongly continuous bundle $\R_+\times\sB(G)$. Note that this does not mean that this strongly continuous bundle is already a strongly continuous product system in the sense Definition \ref{scPSdthm}! For this, the single members $E'_t$ of the bundle (von Neumann \nbd{\cB'}correspondences!) have to embed as right modules into a fixed von Neumann \nbd{\cB'}module $\wh{E'}$, taking the strongly continuous structure from being a subbundle of the trivial strongly continuous bundle $\R_+\times\wh{E'}$, while the $E'_t\subset\sB(G)$ do not sit (for a nontrivial \nbd{E_0}semigroup) in a common von Neumann \nbd{\cB'}submodule contained in $\sB(G)$. (See Section \ref{scsSEC} for details.)

Only running through the theory in \cite[Appendix B]{Ske16} established the intertwiner system, by identifying it thoroughly in Proposition \ref{isoprop} of the present notes algebraically with the strongly continuous commutant of \bf{the} strongly continuous product system of the \nbd{E_0}semi\-group, as a strongly continuous product system in the sense of Definition \ref{scPSdthm}. (The proof of Theorem \ref{ecthm} is, then, done by showing in Lemma \ref{iclem} that the strongly continuous sections of the strongly continuous commutant act, when considered sections of the intertwiner system, strongly continuously on $G$. Before this identification, we would not even have candidates for what the strongly continuous section might be.)

Theorem \ref{ecthm} is, therefore, a result subordinate -- a second order result -- to our theory in \cite{Ske16} of strongly continuous product systems and their strongly continuous commutants. There is no way that Arveson's methods in showing  the same result for \nbd{E_0}semigroups on $\sB(G)$ (and for $t>0$, only) can generalize to the context of \nbd{E_0}semigroups on $\cB\subsetneq\sB(G)$.%
\footnote{ \label{ArvFN}
Arveson, in the proof of \cite[Lemma 2.3]{Arv89}, takes the continuous structure of his intertwiner Arveson system by identifying it with the subspaces $\vt_t(Q)G$ of $G$ where $Q=\gamma\gamma^*$ is a fixed rank-one projection in $\sB(G)$. That is, he identifies his intertwiners with exactly the members of the Bhat system as described in Theorem \ref{Bthm}. While Bhat system and Arveson system of an \nbd{E_0}semigroup on $\sB(G)$ are (algebraically) anti-isomorphic (so that there is no problem to give a continuous structure to the Arveson system by identifying its bundle structure with that of the Bhat system), there is no way to do the same in the general case $\cB\subset\sB(G)$ for the intertwiner system (corresponding to the Arveson system in the case $\cB=\sB(G)$), because its members (\nbd{\cB'}correspondences!) can no longer be identified with the members (\nbd{\cB}correspondences!) of \bf{the} strongly continuous product system (corresponding to the Bhat system in the case $\cB=\sB(G)$).
}
(See also Footnote \ref{AnoncFN}.)

However, now that we know by Theorem \ref{ecthm} that there are enough strongly continuous intertwiner sections for a (faithful) \nbd{E_0}semigroup, we feel encouraged to propose a new, more flexible, definition of \it{strongly continuous product system} in which the intertwiner system (with its strongly continuous structure emerging from being a subbundle of $\R_+\times\sB(G)$) is a strongly continuous product system, too.%
\footnote{
From the beginning with the first definition of \it{continuous product systems} of \nbd{C^*}correspondences (being till now the unchanged basis for all definitions of (strongly) continuous product systems) in Skeide \cite{Ske03b}, we always said that we consider it just a working definition -- a working definition that worked, so far, surprisingly well. (Note that it is not more a working definition than Arveson's measurable product systems, because he even requires that the bundle is isomorphic to a trivial bundle (at the prize of kicking out the point $t=0$ from the bundle; see again Footnote \ref{AnoncFN}), while we only require isomorphism to a subbundle of a trivial bundle.) If the plan of Section \ref{scsSEC} works out, we would finally be prepared to take away the reserve ``working'' from our definition.
}

We give a few hints in Section \ref{scsSEC} how this goes and formulate questions that have to be answered. Also here, a complete treatment has to wait for another paper.
\erem

\brem \label{rem3rd}
We have been asked by a referee of an earlier version (see the acknowledgments) the interesting question whether the situation in these notes could be pushed forward to \it{open quantum systems} (Davies \cite{Dav76}), that is, to \hl{Markov semigroups} (that is, semigroups of (here, normal) unital CP-maps) instead of \nbd{E_0}semigroups.

For a tentative answer to this question, we should say that  Markov semigroups  on $\cB$ are related to \nbd{E_0}semigroups by their so-called (unique) \it{minimal dilation} (an \nbd{E_0}semigroup that can be \it{compressed} to the given Markov semigroup). For the constructions of the dilation, following Bhat and Skeide \cite{BhSk00}, we embed the so-called \it{GNS-correspondences} (Paschke \cite{Pas73}) of the Markov maps into a product system system $E^\odot$ of \nbd{\cB}correspondences $E_t$ with a so-called unit $\xi^\odot$ for $E^\odot$ giving back the Markov semigroup as $\AB{\xi_t,\bullet\xi_t}$. Muhly and Solel  \cite{MuSo02} have taken a similar approach, but starting from the commutant of the GNS-correspondences, the so-called \it{Arveson Stinespring correspondences}, of the Markov maps, also embedding them into a product system, namely, the commutant of $E^\odot$.%
\footnote{
This has been vaguely mentioned in \cite{Ske03c} and several more of our papers; it is clearly stated in Shalit and Skeide \cite[Appendix A(iv)]{ShaSk23}.
}
(The relation between the two ways to obtain the the same \nbd{E_0}semigroup (the minimal dilation is, recall, unique!) is that between the so-called \it{left dilations} of $E^\odot$ and the so-called \it{right dilations} of its commutant; see the discussion of right dilations preceding Theorem and Definition \ref{rdthm} and Footnote \ref{rduniFN}.) Note, however, that there is no unit for the commutant of $E^\odot$ around. (Units for $E^\odot$ go with what \cite{MuSo02} call a \it{completely bounded covariant representation} of the the commutant of $E^\odot$.) So, there is no natural Markov semigroup on $\cB'$ around that may be paired with the orginal one on $\cB$.

Under the extra condition of a cyclic separating vector $\gamma\in G$ such that the vector state $\AB{\gamma,\bullet\gamma}$ on $\cB$ is invariant for the Markov semigroup, there is a \it{dual} Markov semigroup on $\cB'$. (This is a special case of the duality of Markov maps $T\colon\cA\rightarrow\cB$ and $T'\colon\cB'\rightarrow\cA'$ in the presence of covariant states discovered by Albeverio and Hoegh-Krohn \cite{AlHK78}, extended by Gohm \cite{Goh04} to an equivalence of covariant extensions of $T$ to a Markov map $\cA\subset\sB(K)\rightarrow\sB(G)\supset\cB$ and existence of so-called \it{weak tensor dilations} of $T'$. In Gohm and Skeide \cite{GoSk05} we immersed these dualities consistently into the framework of commutants of von Neumann correspondences, and by showing existence of weak tensor dilations, we also showed existence of extensions.)

It is easy to see (using the notations about product systems and their commutants established in our Section \ref{E0psSEC} here) that the GNS-correspondences and the product systems into which they embed are commutants of each other. It makes sense to ask whether the two dilating semigroups of the two Markov semigroups can be paired in some sense. Note, however, that the dilating semigroup of the Markov semigroup on $\cB$ acts on some $\sB^a(E)$ ($E$ a von Neumann \nbd{\cB}module) and the dilating semigroup of the dual Markov semigroup on $\cB'$ acts on some $\sB^a(E')$ ($E'$ a von Neumann \nbd{\cB'}module). Before we can ask whether the two can be paired, we first have to make sit $\sB^a(E)$ and $\sB^a(E')$ in the same $\sB(K)$ as commutants of each other. We believe this is possible. But also this extension of our results here has to wait for another paper.

We should say that the correspondence between Markov semigroups on $\cB$ and on $\cB'$ (with a faithful invariant vector state) is a true duality, while the same \nbd{E_0}semigroup on $\cB$ can be paired with different \nbd{E_0}semigroups on $\cB'$. (By Theorem \ref{cocthm}, two such \nbd{E_0}semigroups on $\cB'$ have to be cocycle equivalent.) So, once the preceding problem of ``embedding'' two dual Markov semigroups and their dilations into one (for positive and negative times) is resolved, we will have to look what this means intrinsically (that is, without running through dilation) for the Markov semigroups, involving possibly also the development of a suitable notion of cocycle equivalence for Markov semigroups. (There are rudimentary considerations abtout the latter in \cite[Section 7]{BhSk00}.) This is still further away than the preceding problem, and has to wait even more.
\erem

\bulletline
We  close this introduction with existing results in special cases and other closely related topics such as extensions of \nbd{E_0}semigroup on $\cB\subset\sB(G)$ to \nbd{E_0}semigroups on $\sB(G)$ that may be automorphic or not. This may safely be skipped.

\brem \label{PRArem}
Powers and Robinson \cite[Definition 3.1]{PoRo89} pair \nbd{E_0}semigroups $\vt_i$ on $\cB_i$ ($i=1,2$) via an automorphism group $\alpha$ on $\cB_1\otimes\cB_2$. Only for $G=H\otimes H$ and $\cB=\sB(H)\otimes\id_H$ this coincides with Definition \ref{paireddefi}. In this situation, Theorem \ref{cpairedthm} amounts to Arveson's result \cite[Theorem 3.5.5]{Arv03} that two strongly continuous \nbd{E_0}semigroups on $\sB(H)$ can be paired if and only if they have anti-isomorphic Arveson systems. Our present notes show, therefore, that Definition \ref{paireddefi} might be better suited to generalize the results by Powers and Robinson and by Arveson from type I factors to general von Neumann algebras.
\erem

A pairing is, in the first place, an extension of $\vt$ to an automorphism group $\alpha$. (This is what Arveson, in the situation of the preceding remark, calls a \it{history}; see \cite[Section 3.5]{Arv03}.) It is noteworthy, that existence of such an extension is already enough to see that such $\vt$ can be paired with \bf{some} $\vt'$. Namely:

\bprop \label{bashisprop}
\hfill
$\alpha_{-t}(\cB)\subset\cB$ ~~~ $\Longleftrightarrow$ ~~~ $\alpha_t(\cB')\subset\cB'$.
\hfill\hfill{~}
\eprop

\proof
This follows since $\alpha_{-t}(b)b'=b'\alpha_{-t}(b)$ for all $b\in\cB,b'\in\cB'$, under $\alpha_t$ transforms into $b\alpha_t(b')=\alpha_t(b')b$ for all $b\in\cB,b'\in\cB'$, and \it{vice versa} under $\alpha_{-t}$.\qed

\lf
So an \nbd{E_0}semigroup $\vt$ on $\cB\subset\sB(G)$ can be paired with \bf{some} \nbd{E_0}semigroup on $\cB'$ if and only if $\vt$ possesses an extension to an automorphism (semi)group on $\sB(G)$. It may possess such an extension or not. (See the examples below.) Of course, for possessing such an extension, the \nbd{E_0}semigroup necessarily has to be faithful. But, whether or not a faithful \nbd{E_0}semigroup does possess an extension, is a question of how the abstract \nbd{W^*}algebra $\cB$ sits as a von Neumann algebra $\cB\subset\sB(G)$. (Indeed, if we consider $\cB$ as an (abstract) \nbd{W^*}algebra (that is, if we do not consider $\cB$ as sitting in $\sB(G)$ from the beginning, so that, in particular, the commutant $\cB'$ is undefined), then by Arveson and Kishimoto \cite{ArKi92} or Skeide \cite[Theorem B.36]{Ske16}, there exists a normal faithful unital representation $\pi\colon\cB\rightarrow\sB(G)$ such that $\vt$ extends from $\cB=\pi(\cB)\subset\sB(G)$ to an automorphism (semi)group on $\sB(G)$, provided $\vt$ is faithful.)

Our theorems are results about (concrete!)\ von Neumann algebras, not about (abstract!)\ \nbd{W^*}algebras.

\bex \label{nonextex}
~\nopagebreak
\begin{enumerate}
\item \label{n1}
No proper \nbd{E_0}semigroup on $\cB=\sB(G)$ can be paired with an \nbd{E_0}semigroup on $\cB'=\id_G\C$. (This follows, for instance, from Theorem \ref{pairedthm}, because, as explained in Section \ref{1dimAsSEC}, an \nbd{E_0}semigroup on $\sB(G)$ is proper if and only if its Arveson system is not one-dimensional, while the Arveson system of the only possible \nbd{E_0}semigroup $\id$ on $\sB(G)'=\id_G\C\cong\sB(\C)$ has as Arveson system the trivial one.)

\item \label{n2}
Take $\cB\subset\sB(G)$ to be a \it{MASA}, so that $\cB'=\cB$. If $\alpha$ is an automorphism of $\sB(G)$ leaving invariant $\cB$, so that by Proposition \ref{bashisprop} also $\alpha^{-1}$ leaves invariant $\cB'=\cB$, then the (co)restriction of $\alpha$ to $\cB$ is invertible. So, no proper endomorphism of $\cB$ can be obtained by (co)restricting an automorphism of $\sB(G)$ to $\cB$, hence no proper \nbd{E_0}semigroup $\vt$ on $\cB$ is the (co)restriction of an automorphism semigroup; that is, no such $\vt$ can be paired. (Concretely, take $\cB:=\C\oplus L^\infty(\R_+)$ acting on $G:=\C\oplus L^2(\R_+)$ and let $\vt$ be the unitalization of the right shift semigroup $\sS$ on $L^\infty(\R_+)$, that is, $\vt_t(z,f):=(z,z\I_{\SB{0,t}}+\sS_t(f))$.) 
\end{enumerate}
\eex

\noindent
It is noteworthy that $\cB=\sB(G)$ in Example \ref{nonextex}\eqref{n1} is not in standard form. (In fact, the \nbd{E_0}semigroups on $\sB(H)$ paired as explained in Remark \ref{PRArem}, live on $G=H\otimes H$ and are, therefore, in standard form.) On the other hand, a MASA is necessarily in standard form; so Example \ref{nonextex}\eqref{n2} yields counter examples in standard form. However, the concrete realization with the unitalized right shift on $\cB=\C\oplus L^\infty(\R_+)$, does not admit \bf{faithful} invariant normal states. (The normal state $(z,f)\mapsto z$ is invariant but, of course, not faithful.) As a MASA, this $\cB$ is the opposite of a factor.  \nbd{E_0}semigroups on factors with invariant normal faithful states, allow statements going much further in the context of extensions to \nbd{E_0}semigroups (not necessarily automorphism semigroups); see, for instance, Bikram, Izumi, Srinivasan, and Sunder \cite{BISS14}. Their setting, where a (proper) \nbd{E_0}semigroup on $\sB(G)$ \bf{also} (co)restricts to an \nbd{E_0}semigroup on $\cB'$ (no $-t<0$ there!)\ is, of course, different from ours here. Still we think that an analysis of their setting in terms of product systems and commutants might be (especially now that we have Theorem \ref{ecthm}!) a promising challenge for the future.

\lf\noindent
\bf{Some basic notations and conventions:}
\it{Homomorphism} between \nbd{*}algebras means \it{\nbd{*}homo\-morphism}. Likewise, for \it{representations}. When there is a product operation $\cA\times\cB\ni(a,b)\mapsto ab\in\cC$ for sets $\cA,\cB,\cC$, for subsets $A\subset\cA$ and $B\subset\cB$ we define the subset
\beqn{
AB
~:=~
\bCB{ab\colon a\in A,b\in B}
}\eeqn
of $\cC$. We \bf{do not adopt} any convention of taking linear spans or even their closure, but will indicate them explicitly when needed.

We assume that the reader knows what a (right, of course) Hilbert module is, while von Neumann modules are explained in Section \ref{E0psSEC}. An element $\xi$ in a Hilbert \nbd{\cB}module is a \hl{unit vector} if $\AB{\xi,\xi}=\U\in\cB$. Recall that a von Neumann module $E$ over a von Neumann algebra $\cB\subset\sB(G)$ is \hl{strongly full} if $\cls^s\AB{E,E}=\cB$, where $\ol{\bullet}^s$ indicates strong closure. We denote the sets of adjointable maps and of bounded and right linear maps on a Hilbert module $E$ by $\sB^a(E)$ and $\sB^r(E)$, respectively. Recall that $\sB^a(E)\subset\sB^r(E)$, in general, and that $\sB^a(E)=\sB^r(E)$ if $E$ is a von Neumann module (or any other self-dual Hilbert module). We use similar notations for maps between Hilbert modules. Especially, for each $x\in E$ we define the map $x^*\in\sB^a(E,\cB)$ by $x^*\colon y\mapsto\AB{x,y}$ (with adjoint map in $\sB^a(\cB,E)$, also denoted by $x$, given by $x\colon b\mapsto xb$). So, $xy^*$ is the \hl{rank-one operator} $z\mapsto x\AB{y,z}$.

Recall that a \hl{\nbd{C^*}correspondence} from $\cA$ to $\cB$ (or an \nbd{\cA}\nbd{\cB}correspondence) is a Hilbert \nbd{\cB}module $E$ that is also an \nbd{\cA}\nbd{\cB}mo\-dule such that the left action acts nondegenerately (that is, $\cls\cA E=E$) and such that the left action defines a (\nbd{*}!)homomorphism (that is $\AB{ax,y}=\AB{x,a^*y}$) into $\sB^a(E)$. For von Neumann correspondences see Section \ref{E0psSEC}. By $\sB^{bil}(E)$ ($\sB^{a,bil}(E)$) we denote the bounded (adjointable) bimodule maps. When $\cA=\cB$, we also will say a correspondence over $\cB$ or a \nbd{\cB}correspondence. The (\hl{internal}) \hl{tensor product} of an \nbd{\cA}\nbd{\cB}cor\-re\-spond\-ence $E$ and a \nbd{\cB}\nbd{\cC}correspondence $F$ is the unique \nbd{\cA}\nbd{\cC}correspondence $E\odot F$ with a total set of elements denoted $x\odot y$ subject to the relations
\baln{
\AB{x\odot y,x'\odot y'}
&
~=~
\AB{y,\AB{x,x'}y'},
&
a(x\odot y)
&
~=~
(ax)\odot y.
}\ealn
For the tensor product of von Neumann correspondences see Section \ref{E0psSEC}. Note that in the particular case when $\cC=\C$ (so that $F$ is a Hilbert space with a nondegenerate representation of $\cB$), there is no difference between the \nbd{C^*}tensor product and the von Neumann tensor product. The same is true, if in a multiple tensor product the last factor is a Hilbert space.

\lf\noindent
\bf{Acknowledgments:} I thank Malte Gerhold, Robin Hillier, and Orr Shalit for useful comments on this and earlier versions, and I thank Malte Gerhold and Michael Schürmann for kind hospitality and support in Greifswald for the last three months of my sabbatical 2019/20, where the algebraic version in Theorem \ref{pairedthm} has been completed.

I also wish to say a big `thank you' to the referee of the present version and to one of the two referees examining an earlier version for CMP. Both have done a superb job.





\newpage

\section{One-dimensional Arveson systems} \label{prelSEC} \label{1dimAsSEC}

As already stated in the introduction, the definition of \bf{the} product system of an \nbd{E_0}semigroup $\vt$ on $\cB$ and the definition of the intertwiner system of $\vt$ (both as given in the introduction) plus referring to Skeide \cite{Ske03c} for that the commutant of the former is (isomorphic to) the latter (also explained in Observation \ref{commob}) plus a \it{few basics} on multipliers and how they act on product systems (also already stated in the introduction), are enough to give a concrete meaning to all pieces that occur in Theorem \ref{pairedthm}. Before proceeding to the proof in Section \ref{pairSEC}, we explain in the present section these \it{few basics} on multipliers, their relation to one-dimensional Arveson systems, and their relation to automorphism (semi)groups on $\sB(G)$.

While surely many statements are \it{folklore} and more or less well-known, the presentation here is quite a bit more than mere repetition. (Most accounts in the literature deal explicitly with measurable multipliers,  which, in the context of strongly continuous automorphism (semi)groups on $\sB(G)$, lead to Wigner's theorem. Non-measurable multipliers, which is precisely why we need this section, are rarely discussed, Liebscher \cite{Lie09} being possibly the only exception; but even there the relations multipliers $\leftrightarrow$ one-dimensional Arveson systems $\leftrightarrow$ automorphism (semi)groups on $\sB(G)$ are not sufficiently explicit for our purposes. Also, as opposed with Arveson's definition (assuming always $t>0$), our definition of product system includes $t=0$ and the related monoid conditions; and this requires some extra considerations not discussed elsewhere.)

\lf
Let us start with recalling Bhat's theorem \cite{Bha96} for Hilbert spaces (the predecessor of the version for Hilbert modules in Theorem \ref{Sthm}, we discuss later) on the construction of an Arveson system for an \nbd{E_0}semigroup on $\sB(G)$ based on selecting a unit vector in $G$.

\bitemp[Theorem \cite{Bha96}.] \label{Bthm}
Let $\Theta$ be an \nbd{E_0}semigroup on $\sB(G)$ and let $\gamma$ be a unit vector in $G$. Then the Hilbert spaces $\eH_t:=\Theta_t(\gamma\gamma^*)G$ with the product maps defined by
\vspace{-1ex}
\beqn{
U_{s,t}
\colon
g_s\otimes h_t
~\longmapsto~
\Theta_t(g_s\gamma^*)h_t
\vspace{-1.5ex}
}\eeqn
form an \hl{Arveson system} (that is, a product system of Hilbert spaces) $\eH^\otimes$, and the  (unitary!) maps defined by
\vspace{-2ex}
\beqn{
V_t
\colon
g\otimes h_t
~\longmapsto~
\Theta_t(g\gamma^*)h_t
}\eeqn
give back $\Theta$ as $\Theta_t=V_t(\bullet\otimes\id_t)V_t^*$. (The $V_t$ form what is called later on a \hl{left dilation} of $\eH^\otimes$ to $G$; see Theorem and Definition \ref{Sthm}.)%
\footnote{
Theorem \ref{Bthm} is one of these results that are fairly easy to check, once one had the (not so easy!) idea to write down the maps $V_t$ and their (co)restrictions $U_{s,t}$. Once their definitions are there, their algebraic properties (isometricity and associativity) are easy to check. And surjectivity is a consequence of normality (so that everything is determined by what it does to rank-one operators).
}
\eitemp

In the situation of Theorem \ref{Bthm}, the members $\Theta_t$ of an \nbd{E_0}semigroup $\Theta$ on $\sB(G)$ are surjective if and only if the members $\eH_t$ of the associated Arveson system $\eH^\otimes$ are one-dimensional. There is an obvious one-dimensional Arveson system, the \hl{trivial} one, with $\eH_t=\C$ and $U_{s,t}(z_s\otimes w_t)=z_sw_t$. Clearly, the Arveson system of the identity semigroup is the trivial one. And since \nbd{E_0}semigroups with the same product system are \it{cocycle equivalent}, and since a cocycle with respect to the identity semigroup is a semigroup itself, we see that an \nbd{E_0}semigroup $\Theta$ has associated the trivial Arveson system if and only if it is induced by a unitary semigroup $U$ as $\Theta_t=U_t\bullet U_t^*$. Such a (semi)group we call \hl{inner}.

By results due to Wigner \cite{Wig39}, all strongly continuous automorphism (semi)groups have this form for a strongly continuous unitary (semi)group $U$; this statement is, therefore, frequently referred to as \hl{Wigner's theorem}. But, there are more complicated one-dimensional Arveson systems. And since every (also non-measurable) Arveson system comes from an \nbd{E_0}semigroup (the proof in Skeide \cite{Ske06} also works without measurability, if we replace direct integrals by direct sums), this means there are (not strongly continuous) automorphism semigroups that cannot be induced by a unitary semigroup (though still by a family of unitaries, but these implementing unitaries need not form a semigroup).

The structure of one-dimensional Arveson systems is determined by so-called \it{multipliers}. Suppose $\eH^\otimes$ is a one-dimensional Arveson system. For each $t$ we choose a unit vector $\gamma_t\in\eH_t$ so that $\eH_t=\gamma_t\C$. A unitary $U_{s,t}\colon\eH_s\otimes\eH_t\rightarrow\eH_{s+t}$ is, then, determined by the unique number $m(s,t)$ in the torus $\T$ such that $U_{s,t}(\gamma_s\otimes\gamma_t)=\gamma_{s+t}m(s,t)$.

\bprop
The function $m\colon\R_+^2\rightarrow\T$ defines a product system structure $U_{s,t}(\gamma_s\otimes\gamma_t)=\gamma_{s+t}m(s,t)$ on the family $\bfam{\gamma_t\C}_{t\in\R_+}$ if and only if
\beq{ \label{mult}
m(r,s)m(r+s,t)
~=~
m(r,s+t)m(s,t),
}\eeq
that is, if $m$ is a \hl{multiplier} (over $\R_+$).
\eprop

\proof
While this condition is, clearly, equivalent to associativity of the product maps, we are left with the marginal conditions, where $s$ or $t$ is $0$. This is dealt with by how we identify $\eH_0=\gamma_0\C$ with $\C=1\C$ as required in the definition of product system (see \ref{Sthm}).

If $U_{0,0}(\gamma_0\otimes\gamma_0)=\gamma_0m(0,0)$, then with $\gamma'_0:=\frac{\gamma_0}{m(0,0)}$ we get $U_{0,0}(\gamma'_0\otimes\gamma'_0)=\gamma'_0$. On the other hand, one easily checks that each multiplier is constant on the set $(\zero\times\R_+)\cup(\R_+\times\zero)$. So, $U_{0,t}(\gamma'_0\otimes\gamma_t)=\gamma_t\frac{m(0,t)}{m(0,0)}=\gamma_t$. Likewise, for $U_{t,0}$. We see, identifying $\gamma'_0\in\eH_0$ with $1\in\C$, the $U_{st}$ define a product system structure.\qed

\bprop
The product system $\eH^\otimes$ is isomorphic to the trivial Arveson system if and only if the multiplier $m$ is \hl{trivial}, that is, if there is a function $f\colon\R_+\rightarrow\T$ such that
\beq{ \label{trimu}
m(s,t)
~=~
\frac{f(s)f(t)}{f(s+t)}.
}\eeq
\eprop

\proof
$\eH$ is isomorphic to the trivial Arveson system if and only if we find unitaries $u_t\colon\gamma_t\C\rightarrow\C,\gamma_t\mapsto f(t)\in\T$ such that
\beqn{
u_{s+t}U_{s,t}(\gamma_s\otimes\gamma_t)
~=~
(u_s\gamma_s)(u_t\gamma_t),
}\eeqn
that is precisely, such that \eqref{trimu} holds.\qed

\brem
It is noteworthy, that the multipliers form a(n abelian) group under pointwise multiplication, and that this group contains $\T$ as the subgroup of constant multipliers. For instance, the multiplier of the $U_{s,t}$ coming by replacing $\gamma_0$ with $\gamma'_0$ is just $m'=\frac{m}{m(0,0)}$. In between sits the subgroup of trivial multipliers. Another subgroup is that consisting of all multipliers $m$ with $m(0,0)=1$. It is isomorphic to the quotient with all constants, and isomorphic to the group of all product system structures on the family $\bfam{\C}_{t\in\R_+}$ under tensor product of  Arveson systems (based on $\eH_t\otimes\eH'_t\ni 1\otimes 1\mapsto 1\in\eH''_t$, where $\eH''^\otimes$ is the product of $\eH^\otimes$ and $\eH'^\otimes$). The quotient of all normalized multipliers by its trivial ones is precisely the group of all one-dimensional Arveson systems up to isomorphism. Note that with $m$ also the transpose $m^\dagger\colon(s,t)\mapsto m(t,s)$ is a multiplier. Clearly, $m\mapsto m^\dagger$ is a self-inverse group automorphism.
\erem

\brem
Most about multipliers can be found in Arveson \cite[Section 3.5]{Arv03}, provided we require the multipliers measurable. But, it is a theorem that measurable multipliers are trivial. Liebscher proves this result in two (interesting!)\ stages: Symmetric multipliers (that is, $m=m^\dagger$) are trivial; measurable multipliers are symmetric. Of course, there are trivial multipliers that are non-measurable. In \cite[Example 7.17]{Lie09}, Liebscher presents a (necessarily non-measurable) multiplier that is not symmetric, hence, not trivial. \cite{Lie09} discusses multipliers over $\R$, not $\R_+$, though, but as the latter always extend as multipliers to $\R^2$, this is not a problem. (This follows, again, from existence of an \nbd{E_0}semigroup for any Arveson system, and, as noted already by Arveson, that for a one-dimensional Arveson system, this \nbd{E_0}semigroup consists of automorphisms and, thus, extends to an automorphism group. Repeating the construction in Theorem \ref{Bthm} for all $t\in\R$, we get a one-dimensional product system over $\R$, and from that a multiplier over $\R$ as before for $\R_+$, which, choosing for $t\ge0$ the same $\gamma_t$ as before, extends the given multiplier over $\R_+$.)
\erem

\lf
We said already in the the introduction that a multiplier $m$ (that is, a one-dimensional Arveson system $\eH^\otimes$) can act on a general product system $E^\odot$, replacing its product $u_{s,t}$ by $x_s\odot y_t\mapsto u_{s,t}(x_s\odot y_t)m(s,t)$. (This is just either of the two external tensor products $E^\odot\otimes\eH^\otimes$ or $\eH^\otimes\otimes E^\odot$. If $E^\odot$ has a trivial subsystem (that is, if $E^\odot$ possesses a so-called \it{central unital unit} or, equivalently, if $E^\odot$ is \it{spatial}; see Skeide \cite{Ske06d}), then $m$ may be recovered, so the group of multipliers acts faithfully on $E^\odot$. If $E^\odot$ is nonspatial, we do not know.)

The formulation of Theorem \ref{pairedthm} means that there is a one-dimensional Arveson system $\eH^\otimes$ (or a multiplier) such that the commutant of the product system of $\vt$ is isomorphic to the product system of $\vt'$ tensored with $\eH^\otimes$.

The following result is most probably \it{folklore}. But we do not know a reference.

\bthm \label{nonWigthm}
If $\alpha$ is a (semi)group of automorphisms of $\sB(G)$ (automatically, normal), then there is a family of unitaries $U_t\in\sB(G)$ and a multiplier $m$ on $\R$ (on $\R_+$) such that
\beqn{
U_tU_s
~=~
U_{s+t}m(s,t)
\vspace{-2ex}
}\eeqn
and such that $\alpha_t=U_t\bullet U_t^*$.
\ethm

\proof
Do the construction in Theorem \ref{Bthm} for $\Theta=\alpha$ (if necessary extended to $t\in\R$) and produce $\eH_t=\gamma_t\C$, $U_{s,t}$, and $V_t$ (if necessary for all $s,t\in\R$). Put $U_t\colon g\mapsto V_t(g\otimes\gamma_t)$. Then $U_t\bullet U_t^*=\alpha_t$ and
\beqn{
U_tU_sg
~=~
V_t\Bfam{V_s(g\otimes\gamma_s)\otimes\gamma_t}
~=~
V_{s+t}\Bfam{g\otimes U_{s,t}(\gamma_s\otimes\gamma_t)}
~=~
V_{s+t}(g\otimes\gamma_{s,t}m(s,t))
~=~
U_{s+t}gm(s,t).\qedsymbol
}\eeqn
\noqed

\vspace{-3ex}
\brem \label{opprem}
Note that $U_tU_s$ goes with  $m(s,t)$, that is, with $m^\dagger$, not with $m$. (As multipliers need not be symmetric, this is not marginal.) This is okay so, in several ways. First, as noticed in Shalit and Skeide \cite[Section 4]{ShaSk23}, the product system \it{à la} Bhat goes, for non-abelian indexing monoids, with the opposite monoid of the semigroup from which it is derived. Second, the Bhat system of an \nbd{E_0}semigroup is the opposite of the system constructed by Arveson (see Skeide \cite[Section 2]{Ske05a}), because the latter is actually the commutant of the former, and the commutant is anti-multiplicative for the tensor product (see Section \ref{E0psSEC} and Skeide \cite[Example 6.4]{Ske08a}). Third, the careful reader will notice that, in the following section, the product system goes with $\vt$, and in a pairing $\vt_t$ goes, recall!, with $\alpha_{-t}$. Therefore, in the following section we will have $U_sU_t$ going with $m(s,t)$.
\erem

\newpage
\section{Proof of Theorem \ref{pairedthm}} \label{pairSEC}

Recall that, in the situation of Theorem \ref{pairedthm}, \bf{the} product system of $\vt$ on $\cB$ and \bf{the} product system of $\vt'$ on $\cB'$ are commutants of each other (modulo a multiplier) if one is isomorphic (modulo a multiplier) to the commutant of the other. We shall examine \bf{the} product system $E^\odot=\bfam{E_t}_{t\in\R_+}$ of $\vt$ and the commutant $F^\odot=\bfam{F_t}_{t\in\R_+}$ of \bf{the} product system of $\vt'$. Recall that the former is given as $E_t=\cB_t={_{\vt_t}}\cB$ and comes along with the \it{left dilation} $v_t\colon b\odot x_t\mapsto\vt_t(b)x_t=b.x_t$ to $E=\cB$ giving back $\vt$ as $\vt_t=v_t(\bullet\odot\id_t)v_t^*$. And recall that the latter comes along as the concrete product system of von Neumann \nbd{\cB}correspondences
\beq{ \label{Bintdef}
F_t
~=~
\bCB{y_t\in\sB(G)\colon\vt'_t(b')y_t=y_tb'~(b'\in\cB')}
~\subset~
\sB(G)
}\eeq
in the identity representation. There is also the \it{right dilation} $w_t\colon y_t\odot g\mapsto y_tg$ to $H=G$ giving back $\vt'$ as $\vt'_t=w_t(\id_t\odot\bullet)w_t^*$.%
\footnote{ \label{lrdFN}
The formal definition of  \it{left dilation} (and also of \it{product system}) is stated in Theorem and Definition \ref{Sthm}, while that of \it{right dilation} is stated in Theorem and Definition \ref{rdthm}. For the following proof we only need to know that the (obviously isometric and left linear) $v_t$ and $w_t$ are unitary (obvious for $v_t$, and well known for $w_t$ as intertwiner spaces act nondegenerately), and the associativity conditions $v_t(v_s\odot\id_t)=v_{s+t}(\id_{E=\cB}\odot u_{s,t})$ (easy) and $w_s(\id_s\odot w_t)=w_{t+s}(v_{s,t}\odot\id_{H=G})$ (obvious) required in Equation \eqref{Usgmeq}.
}

\lf
\proof[Proof of Theorem \ref{pairedthm} ''$\Leftarrow$''.~]
Now, if $E^\odot$ and $F^\odot$ are isomorphic modulo a multiplier $m$ (which we assume, for simplicity, normalized as $m(0,0)=1$), this means there are isomorphisms $u_t\colon E_t\rightarrow F_t$ such that $u_0=\id_\cB\colon E_0=\cB\rightarrow\cB=F_0$ and
\beq{ \label{misodef}
u_{s+t}u_{s,t}(x_s\odot\hat{x}_t)=v_{s,t}(u_sx_s\odot u_t\hat{x}_t)m(s,t),
}\eeq
where $u_{s,t}$ (recall, $x_s\odot\hat{x}_t\mapsto\vt_t(x_s)\hat{x}_t$) and $v_{s,t}$ (recall, $y_s\odot\hat{y}_t\mapsto y_s\hat{y}_t$) denote the product maps of $E^\odot$ and $F^\odot$\!, respectively. Using $u_t$, $v_t$ and $w_t$, on $\cB\odot G=G$ we may define the unitaries
\beq{ \label{Utdef}
U_t
\colon
v_t(x\odot x_t)\odot h
~\longmapsto~
x\odot w_t(u_tx_t\odot h)
}\eeq
for $x\in\cB$, $x_t\in E_t$, $h\in G$. We show that $U_{s+t}=U_sU_tm(s,t)$. Indeed,
\bmu{ \label{Usgmeq}
\bfam{U_sU_tm(s,t)}U_{s+t}^*\Bfam{x\odot w_{s+t}\Bfam{v_{s,t}(u_sy_s\odot u_tz_t)\odot h}}
\\
~=~
U_sU_tU_{s+t}^*\Bfam{x\odot w_{s+t}\Bfam{u_{s+t}u_{s,t}(y_s\odot z_t)\odot h}}
~=~
U_sU_t\Bfam{v_{s+t}\Bfam{x\odot u_{s,t}(y_s\odot z_t)}\odot h}~~~~~\,
\\
{~}\hspace{3.5ex}
~=~
U_sU_t\Bfam{v_t\Bfam{v_s(x\odot y_s)\odot z_t}\odot h}
~=~
U_s\Bfam{v_s(x\odot y_s)\odot w_t(u_tz_t\odot h)}
\\
~=~
x\odot w_s\Bfam{u_s y_s\odot w_t(u_tz_t\odot h)}
~=~
x\odot w_{s+t}\Bfam{v_{s,t}(u_s y_s\odot u_tz_t)\odot h},
}\emu
where in the passage from the first  to the second line we used \eqref{misodef}, and in the passage from the second to the third line and in the last equality we used the associativity properties mentioned in Footnote \ref{lrdFN}.%
\footnote{
The definition of the $U_t$ and the verification that they form a semigroup modulo $m$ works for any left dilation $v_t$ of $E^\odot$ to $E$ and any right dilation $w_t$ of $F^\odot$ to $H$. (This is the reason why also in our special situation where $E=\cB$ and $H=G$ we preferred to write letters $x\in \cB$ and $h\in G$.) This generalizes the construction in our proof of the Arveson-Kishimoto theorem \cite{ArKi92} (which corresponds  to $F^\odot=E^\odot$; see \cite[Theorem 3.1]{Ske11a} or \cite[Theorem B.36]{Ske16}) to the situation here where $F^\odot$ is only isomorphic via a multiplier to $E^\odot$.
}
So $\alpha_t:=U_t\bullet U_t^*$ defines a semigroup. Directly from \eqref{Utdef}, applying $\alpha_t(\id_\cB\odot b')$ to $x\odot w_t(u_tx_t\odot h)$ and $\alpha_{-t}(b\odot\id_G)$ to $v_t(x\odot x_t)\odot h$ ($t\ge0$), we see that $\alpha$ pairs $\vt$ and $\vt'$.\qed

\lf
\proof[Proof of Theorem \ref{pairedthm} ''$\Rightarrow$''.~]
Suppose that $\alpha$ is a pairing for $\vt$ and $\vt'$. We know from Theorem \ref{nonWigthm} that $\alpha$ has the form $\alpha_t=U_t\bullet U_t^*$ for unitaries $U_t\in\sB(G)$ that form a semigroup modulo a multiplier $m$. To be compatible with the preceding part, we choose the multiplier such that $U_{s+t}=U_sU_tm(s,t)$. (See also Remark \ref{opprem}.)

The intertwiner condition that characterizes the elements of $F_t$, in terms of $U_t$ reads
\beqn{
U_tb'U_t^*y_t
~=~
y_tb'
}\eeqn
or, equivalently, $b'U_t^*y_t=U_t^*y_tb'$ for all $b'\in\cB'$. Therefore, $y_t\in F_t$ if and only if $U_t^*y_t\in\cB$, so $F_t=U_t\cB$. Recall that the left action of $\cB$ on $F_t$ is $\id_\cB$. So, for $y_t=U_tb_1\in F_t$ we find $b_2y_t=b_2U_tb_1=U_t\vt_t(b_2)b_1=U_t(b_2.b_1)$. It follows that $u_t\colon x_t\mapsto U_t x_t$ is an isomorphism $E_t\rightarrow F_t$ of von Neumann \nbd{\cB}correspondences.

By
\vspace{-2ex}
\bmun{
\hspace{7ex}
u_{s+t}\Bfam{u_{s,t}(x_s\odot \hat{x}_t)}
~=~
U_{s+t}\vt_t(x_s)\hat{x}_t
~=~
U_sU_t\vt_t(x_s)\hat{x}_tm(s,t)
\\
~=~
U_sx_sU_t\hat{x}_tm(s,t)
~=~
v_{s,t}(u_sx_s\odot u_t\hat{x}_t)m(s,t)
\hspace{7ex}
}\emun
we see that the $u_t$ form an isomorphism modulo $m$ of product systems.\qed

\lf
This completes the proof of Theorem \ref{pairedthm}.

\lf
Note that, given $\vt$ and $\vt'$, in the proof we have established a one-to-one correspondence between unitary (semi)groups modulo a multiplier $U_t$ in $\sB(G)$ such that $\alpha_t:=U_t\bullet U_t^*$ is a pairing of $\vt$ and $\vt'$, and isomorphisms modulo the same multiplier $u_t$ between the product system $E^\odot$ of $\vt$ and the commutant $F^\odot$ of the product system of $\vt'$, via Equation \ref{Utdef}. (The inverse is given by $u_t=U_t$, meaning that concrete operator multiplication by $U_t$ defines a bilinear unitary $u_t$.) Since the multiplier on each side is the same, this means, in particular, that if it is absent ($=$ constant $1$) on one side then so it is on the other.

Calling a pairing \hl{inner} if the $U_t$ implementing it can be chosen to form a group, we obtain the following intermediate result (lying strictly in between Theorem \ref{pairedthm} and Theorem \ref{cpairedthm}):

\bthm \label{ipairedthm}
An \nbd{E_0}semigroup $\vt$ on $\cB$ and an \nbd{E_0}semigroup $\vt'$ on $\cB'$ can be paired via an inner pairing $\alpha$ if and only if their product systems are commutants of each other.
\ethm

This is an important step towards the proof of Theorem \ref{cpairedthm} in Section \ref{cpairSEC}.

\section{Strongly continuous product systems: Proof of Theorem \ref{ecthm}} \label{E0psSEC}

Recall that to prove Theorem \ref{cpairedthm}, we need Theorem \ref{pairedthm} (especially, in the particular case in Theorem \ref{ipairedthm}) and Theorem \ref{ecthm}. So far, to understand and prove Theorem \ref{pairedthm} we could content ourselves with the definition of \bf{the} product system and the intertwiner product system of an \nbd{E_0}semigroup on $\cB$ as stated in the introduction, and with little knowledge from Skeide \cite{Ske03c} that states that the commutant of the former is (isomorphic to) the latter. At that stage, we are also able to interpret Theorem \ref{ecthm} (making just a statement about the sections of the intertwiner system). However, we did not yet introduce all notions occurring in Theorem \ref{cpairedthm}. In this section, we provide the notions of \it{strongly continuous product systems} and their \it{strongly continuous commutants} and several results about them from Skeide \cite{Ske16}, necessary to interpret Theorem \ref{cpairedthm} and to prove Theorems \ref{ecthm} and \ref{cpairedthm}. The section culminates in the proof of Theorem \ref{ecthm} (while the proof of Theorem \ref{cpairedthm} has to wait for the following section). For that proof, we need to understand very explicitly the (algebraic) isomorphism between  the strongly continuous commutant of \bf{the} product system and the intertwiner system, which has never been done before.

Before we say what a strongly continuous product system is and what \bf{the} strongly continuous product system of a strongly continuous \nbd{E_0}semigroup is (\cite[Chapter 12]{Ske16}), let us outline the steps how to construct from it (provided product system and \nbd{E_0}semigroup are faithful) its strongly continuous commutant (\cite[Appendix B]{Ske16}):

\bemp[Outline of the construction of the strongly continuous commutant.~] \label{outob}
\begin{itemize}
\item
After having understood the strongly continuous structure of \bf{the} product system $E^\odot$ of a strongly continuous \nbd{E_0}semigroup $\vt$ (\cite[Section 12]{Ske16}), provided $\vt$ is faithful (a minimum requirement for being pairable) construct a strongly continuous so-called \it{right dilation} of $E^\odot$ (granted by \cite[Theorem B.31]{Ske16}).

\item
The right dilation comes along with a normal nondegenerate representation $\rho$ of $\cB$ on a Hilbert space $H$ and it induces a strongly continuous \nbd{E_0}semigroup $\vt'$ acting on $\rho(\cB)'\subset\sB(H)$. This von Neumann algebra $\rho(\cB)'$ is nothing but the algebra $\sB^a(E')$ for the (strongly full) \it{concrete von Neumann \nbd{\cB'}module} $E':=C_\cB(\sB(G,H))$. (Recall Footnote \ref{CBMdefFN}.)

\item
The \it{strongly continuous commutant} of $E^\odot$, $F'^\odot$, is, then, defined as \bf{the} strongly continuous product system of the \nbd{E_0}semigroup $\vt'$ on $\sB^a(E')$.

(Note that we called the strongly continuous commutant $F'^\odot$ in order to not confuse it with the algebraic commutant of the product system (of concrete von Neumann correspondences) $E^\odot$, $E'^\odot$, which is defined in a different way as intertwiner system. We repeat: A good deal of the new work to be done in these notes consists in writing down carefully an explicit (algebraic) isomorphism.)
\end{itemize}
\eemp

\noindent
Of course, this description is more than rushy and leaves several things to be explained. (\bf{The} product system of an \nbd{E_0}semigroup $\vt$ on $\sB^a(E)$ provided $E$ has a unit vector $\xi$; its strongly continuous structure when $\vt$ is strongly continuous; a precise definition of the algebraic commutant; and for the latter two, we need to be specific about what we mean by (concrete) von Neumann modules and correspondences.) The reason we outlined the construction of the strongly continuous commutant already here, is the insight that in order to understand the strongly continuous commutant we need to understand the construction of the strongly continuous structure of \bf{the} product system of a strongly continuous \nbd{E_0}semigroup not only in the particularly simple case where $\vt$ acts on $\cB$ itself (where $\cB=\sB^a(\cB)$ for the von Neumann \nbd{\cB}module $\cB$), but in full generality (where $\vt$ may act on $\sB^a(E)$ for a strongly full von Neumann \nbd{\cB}module $E$), although the construction will, then, be applied to $\vt'$ on $\sB^a(E')$ rather than to $\vt$ itself: We need the full theory from \cite[Section 12+Appendix B]{Ske16} -- and, for Theorem \ref{ecthm}, we need to go beyond \cite{Ske16}.

\bulletline[The product system of an \nbd{E_0}semigroup.~]
The construction of a product system for an \nbd{E_0}semigroup on $\sB^a(E)$ based on the assumption to have a unit vector $\xi\in E$, goes back to Bhat's construction \cite{Bha96} for $E$ being a Hilbert space (we discussed already in Theorem \ref{Bthm}) and its generalization to Hilbert modules $E$ in Skeide \cite{Ske02}. Once a unit vector is granted, the construction does not really depend on the category, provided the endomorphisms $\vt_t$ in the \nbd{E_0}semigroup are sufficiently regular and the tensor product is that within the category under consideration: For Hilbert modules, this condition is that each $\vt_t$ be strict; for von Neumann or \nbd{W^*}modules the $\vt_t$ have to be only normal.%
\footnote{
In Muhly, Skeide, and Solel \cite{MSS06} it has been explained that these conditions amount to a nondegeneracy condition for the action of $\vt_t(\sF(E))$ ($\sF(E)$ denoting the finite-rank operators) on $E$. (For $\vt_t$ being strict, we need that $\vt_t(\sF(E))E$ is norm-total in $E$; for normality it suffices that this set is total in any of the (\nbd{\sigma})weak or (\nbd{\sigma})strong topologies of the von Neumann or \nbd{W^*}module $E$.)
}
Whatever version the reader knows among (concrete) von Neumann or \nbd{W^*}modules and correspondences and their tensor products, the following theorem about \nbd{E_0}semigroups (by definition, recall!, normal) on the von Neumann or \nbd{W^*}algebra $\sB^a(E)$ is a simple adaptation of the \nbd{C^*}result in \cite{Ske02}:

\bitemp[Theorem and Definition \cite{Ske02}.~] \label{Sthm}
Let $E$ be a (concrete) von Neumann or \nbd{W^*}module over a von Neumann or \nbd{W^*}algebra $\cB$ with a unit vector $\xi\in E$, and let $\vt$ be an \nbd{E_0}semigroup on $\sB^a(E)$. Then the von Neumann or \nbd{W^*}correspondences $E_t:=\vt_t(\xi\xi^*)E$ over $\cB$ with left action defined by
\beqn{
b.x_t
~:=~
\vt_t(\xi b\xi^*)x_t
}\eeqn
and product maps determined by
\beqn{
u_{s,t}
\colon
x_s\odot y_t
~\longmapsto~
\vt_t(x_s\xi^*)y_t
}\eeqn
form a \hl{product system} $E^\odot=\bfam{E_t}_{t\in\R_+}$, that is: The $u_{s,t}$ are isomorphisms $E_s\sodots E_t\rightarrow E_{s+t}$ of correspondences (in the relevant category with the tensor product $\sodots$ in that category), the product $(x_s,y_t)\mapsto x_sy_t:=u_{s,t}(x_s\odot y_t)$ is associative, and the marginal conditions $E_0=\cB$ and $u_{0,t}$ and $u_{t,0}$ being the canonical identifications $\cB\odot E_t\equiv E_t\equiv E_t\odot\cB$ are fulfilled. We call $E^\odot$ {\normalfont\bf{the}} \hl{product system} of $\vt$ associated with $\xi$.

Moreover, the maps $v_t\colon E\sodots E_t\rightarrow E$ determined by
\beqn{
v_t
\colon
x\odot y_t
~\longmapsto~
\vt_t(x\xi^*)y_t
}\eeqn
form what is called a \hl{left dilation} of $E^\odot$ to $E$ (namely, $E$ is strongly full, the $v_t$ are unitary, and the product $xy_t:=v_t(x\odot y_t)$ iterates associatively with the product of the product system) giving back $\vt$ as $\vt_t=v_t(\bullet\odot\id_t)v_t^*$.
\eitemp

\bob \label{TPSob}
We note the following.
\begin{itemize}
\item
Obviously, for $E=\cB$ and $\xi=\U\in\cB$ we recover \bf{the} product system of $\vt$ on $\cB$ as described in the introduction (in the second item of the discussion following Theorem \ref{pairedthm}).

\item
Whenever for a given product system $E^\odot$ we have a left dilation $v_t$ to a (by definition) strongly full $E$, then by setting
\beqn{
\vt^v_t
~:=~
v_t(\bullet\odot\id_t)v_t^*
}\eeqn
we define an \nbd{E_0}semigroup $\vt^v$ on $\sB^a(E)$. Note that for admiting such a left dilation, the members $E_t$ of $E^\odot$ have to be necessarily strongly full.

\item
If $E$ has a unit vector, then \bf{the} product system of $\vt^v$ formed by $\vt^v_t(\xi\xi^*)E$ is canonically isomorphic to $E^\odot$ via $ E_t\ni x_t\mapsto \xi x_t=v_t(\xi\odot x_t)\in\vt^v_t(\xi\xi^*)E$.

We see:  Every left dilation $v_t$ of a product system $E^\odot$ to a given $E\ni\xi$ that gives back a given \nbd{E_0}semigroup $\vt$ on $\sB^a(E)$ as $\vt=\vt^v$, identifies (up to isomorphism) $E^\odot$ as \bf{the} product system of $\vt$.%
\footnote{ \label{PSuniFN}
This is even true without unit vector $\xi$ as long as $E$ is (as required for a left dilation) strongly full; see \cite[Proposition 6.3]{Ske16}. In these notes we omit the general case, as we do not need it.
}
\end{itemize}
\eob
The last fact in the preceding observation is responsible for why to quite an extent we needed not worry in our earlier papers about how the product system of an \nbd{E_0}semigroup is actually constructed; it is enough that we also construct a left dilation that gives back the \nbd{E_0}semigroup we started with. But (as first observed in Skeide \cite{Ske03b} for the \nbd{C^*}case), as soon as strong continuity is concerned, it is the construction of Theorem \ref{Sthm} which we have to use and which allows to indicate the strongly continuous structure of \bf{the} product system of an \nbd{E_0}semigroup. This is so, simply because the construction based on a unit vector $\xi$ comes along with an identification of $E_t$ as the subset $\vt_t(\xi\xi^*)E$ of $E$, allowing us to take the continuity of sections as mappings into $E$. Already on the level of the \nbd{E_0}semigroup on $\sB^a(E)$, we have to say what the strong topology of $\sB^a(E)$ is, when $\cB$ is not given \it{a priori} as a von Neumann algebra (that is, if $E$ is a von Neumann module), but as a \nbd{W^*}algebra (when $E$ is only a \nbd{W^*}module); and it is definitely the \bf{strong} topology of a von Neumann module that allows (thanks to joint strong continuity on bounded subsets of operator multiplication) to get strong continuity of the product of \bf{the} product system of a strongly continuous \nbd{E_0}semigroup easily and in the most painless and natural way:%
\footnote{
One may take the \nbd{\sigma}strong topology, that is, the strong topology in a ``big enough'' representation, where ``big enough'' means that every normal state is a vector state. On the level of \nbd{W^*}modules this is the \nbd{s}topology introduced by Baillet, Denizeau, and Havet \cite{BDH88}. But it is cumbersome, and not really efficient to understand and exploit the interplay between these  ``\nbd{W^*}strong topologies'' and how continuity can be check in them.
}

\bulletline[The strongly continuous product system of a strongly continuous \nbd{E_0}semigroup.~]
For dealing with the strongly continuous structure, we continue with some basics about von Neumann modules. (There is a quick introduction to von Neumann modules and von Neumann correspondences in Shalit and Skeide \cite[Appendix A]{ShaSk23}.)
\begin{itemize}
\item
A \hl{concrete} (pre-)Hilbert module over a concrete (pre-)\nbd{C^*}algebra $\cB\subset\sB(G)$ ($G$ a Hilbert space) is a (linear) subspace $E\subset\sB(G,H)$ ($H$ another Hilbert space) such that
\baln{
\cls EG
&
~=~
H,
&
E\cB
&
~\subset~
E,
&
E^*E
~\subset~
\cB.
}\ealn
By the second and the third condition, $E$ with inner product $\AB{x,y}:=x^*y$ and right multiplication $(x,b)\mapsto xb$ (operator multiplication) is a (pre-)Hilbert \nbd{\cB}module. Although this is not strictly necessary, we usually will assume that $\cB$ acts nondegenerately on $G$.

\item
Every (pre-)Hilbert module $E$ over a concrete (pre-)\nbd{C^*}algebra $\cB\subset\sB(G)$ may be turned into a concrete one: For each $x\in E$ define $L_x\in\sB(G,E\odot G)$ by $L_xg:=x\odot g$. Then $x\mapsto L_x$ is an isomorphism from $E$ onto the concrete (pre-)Hilbert \nbd{\cB}module
\beqn{
\bCB{L_x\colon x\in E}
~\subset~
\sB(G,E\odot G).
}\eeqn
Moreover, if $u\colon E\rightarrow F$ is another isomorphism onto a concrete (pre-)Hilbert \nbd{\cB}module $F\subset\sB(G,H)$, then $U\colon x\odot g\mapsto (ux)g$ is a unitary $E\odot G\rightarrow H$ (satisfying $UL_x=ux$, and determined uniquely by that).%
\footnote{
$U$ is, clearly, an isometry. For surjectivity we see the nondegeneracy condition $\cls FG=H$ at work.
}

Clearly, $\sB^a(E)$ sits as $\sB^a(E)\odot\id_G$ in $\sB(E\odot G)$, and so does $\sB^a(F)$ in $\sB(H)$.

\item
A concrete (pre-)Hilbert module $E\subset\sB(G,H)$ over a von Neumann algebra $\cB\subset\sB(G)$ is a \hl{concrete von Neumann \nbd{\cB}module} if $E$ is strongly closed in $\sB(G,H)$; see Skeide \cite[Definition 3.9]{Ske06b}. We denote this as $(E,H)$.

A (pre-)Hilbert module $E$ over a von Neumann algebra $\cB\subset\sB(G)$ is a \hl{von Neumann \nbd{\cB}module} if its ``concretization'' as $E\subset\sB(G,E\odot G)$ is a concrete von Neumann module or, equivalently, if the \hl{extended linking algebra} $\rtMatrix{\cB&E^*\\E&\sB^a(E)}=\sB^a\rtMatrix{\cB\\E }\subset\sB\rtMatrix{G\\E\odot G}$ is a von Neumann algebra on $\rtMatrix{G\\E\odot G}$; see Skeide \cite[Definition 4.4 and Proposition 4.5]{Ske00b}.%
\footnote{
Since a self-dual (pre-)Hilbert module $E$ over a von Neumann algebra $\cB$ is necessarily strongly closed in $\sB(G,E\odot G)$, by \cite[Theorem 4.16]{Ske00b} we see that a (pre-)Hilbert \nbd{\cB}module is a von Neumann \nbd{\cB}module if and only if it is self-dual, that is, if and only if it is a \nbd{W^*}module. Of course, after identifying a \nbd{W^*}algebra $\cB$ as a von Neumann algebra on $G$ (via a nondegenerate faithful representation $\pi\colon\cB\rightarrow\sB(G)$), every \nbd{W^*}module over $\cB$  is a (unique up to suitable unitary equivalence, as soon as $\pi$ is fixed) von Neumann module over $\cB=\pi(\cB)\subset\sB(G)$.
}
\end{itemize}

\noindent
We are now ready to discuss \bf{the} strongly continuous product system of a strongly continuous \nbd{E_0}semigroup. (We formulate it for von Neumann modules and correspondences. In this form, it applies, of course, also to concrete ones, but for concrete ones there is an obvious unitarily equivalent version where $E_t\odot G$ and $\wh{E}\odot G$ would be replaced by $H_t$ (for $(E_t,H_t)$) and $\wh{H}$ for $(\wh{E},\wh{H}$), respectively.) We again leave to reader's discretion from which category -- von Neumann or concrete von Neumann -- modules, correspondences, and tensor products are. Only when we come to commutants, we will fix on the \it{concrete} category and its tensor product is explained later on.

\bitemp[Definition and Theorem {\cite[Chapter 12]{Ske16}}.~] \label{scPSdthm}
Let $\cB\subset\sB(G)$ be a von Neumann algebra on $G$ and let $E^\odot$ be a product system of von Neumann \nbd{\cB}correspondences. For a von Neumann \nbd{\cB}module $\wh{E}$ and a family $i_t\colon E_t\rightarrow\wh{E}$ of isometric embeddings (as right module; $\wh{E}$ need not be a correspondence), we denote by
\beqn{
CS_i^s(E^\odot)
~:=~
\BCB{\bfam{x_t}_{t\in\R_+}\colon x_t\in E_t, t\mapsto i_tx_t\odot g\in\wh{H}:=\wh{E}\odot G\text{\normalfont~is continuous for all~}g\in G}
}\eeqn
the (right \nbd{\cB}module of) \hl{strongly continuous sections} of $E^\odot$ (with respect to $i_t$). We say $E^\odot$ is a \hl{strongly continuous product system} (with strongly continuous structure $i_t$) if
\beqn{
\BCB{x_s\colon\bfam{x_t}_{t\in\R_+}\in CS_i^s(E^\odot)}
~=~
E_s
}\eeqn
for all $s\in\R_+$ (that is, if $E^\odot$ has \hl{enough strongly continuous sections}), and if the function
\beqn{
(s,t)
~\longmapsto~
i_{s+t}(x_sy_t)\odot g\in\wh{H}
}\eeqn
is continuous for all $\bfam{x_t}_{t\in\R_+},\bfam{y_t}_{t\in\R_+}\in CS_i^s(E^\odot)$ and for all $g\in G$.

A morphism between strongly continuous product systems is \hl{strongly continuous} if it sends strongly continuous sections to strongly continuous sections. By \cite[Theorem B.5]{Ske16}, a strongly continuous isomorphism, has a strongly continuous inverse.

{\normalfont\bf{The}} product system $E^\odot$ of a strongly continuous \nbd{E_0}semigroup $\vt$ on $\sB^a(E)$ for a von Neumann \nbd{\cB}module $E$ with a unit vector $\xi$ and with $i_t$ the canonical embedding of $E_t=\vt_t(\xi\xi^*)E\subset E$ into $E$ is strongly continuous. We call $E^\odot$ with this strongly continuous structure {\normalfont\bf{the}} \hl{strongly continuous product system} of $\vt$.

The left dilation $v_t$ is \hl{strongly continuous} in the sense that $t\mapsto(xy_t)\odot g\in E\odot G$ is continuous for all $x\in E$, all $\bfam{y_t}_{t\in\R_+}\in CS_i^s(E^\odot)$ and all $g\in G$. Moreover (see \cite[Observation B.14]{Ske16}), if $v_t$ is any strongly continuous left dilation of any strongly continuous product system $E^\odot$ to any $E\ni\xi$, then the \nbd{E_0}semigroup $\vt^v$ is strongly continuous and {\normalfont\bf{the}} strongly continuous product system of $\vt^v$ is $E^\odot$.
\eitemp

In the item before their definition, we have seen that there is not much of a difference between von Neumann modules $E\subset\sB(G,E\odot G)$ and concrete von Neumann modules $E\subset\sB(G,H)$. However, while it is very important to never forget that $H$ actually \bf{is} (canonically isomorphic to) the tensor product $E\odot G$ (with all its algebraic implications which, otherwise, ``fall from heaven''), the picture as concrete von Neumann modules has the striking advantage that certain functors become one-to-one functors (not just natural equivalences).

An instance of this ambivalence leads straight to the definition of commutant:

\bulletline[The algebraic commutant.~]
On $E\odot G$ we have the action of $\sB^{bil}(G)=\cB'$ as $b'\mapsto\id_E\odot b'$, which in the concrete picture $E\subset\sB(G,H)$ acts as $\rho'(b')\colon xg\mapsto xb'g$. We refer to the (normal and nondegenerate) representation $\rho'$ as the \hl{commutant lifting} of $\cB'$ associated with $E\subset\sB(G,H)$. Conversely, if $\rho'$ is a (normal and nondegenerate) representation of $\cB'$ on $H$, then $E:=C_{\cB'}(\sB(G,H))$ is a concrete von Neumann \nbd{\cB}module.%
\footnote{
The only not entirely obvious issue, nondegeneracy $\cls EG=H$, follows from normality of $\rho'$. (While this is surely \it{folklore}, possible references are Muhly and Solel \cite[Lemma 2.10]{MuSo02} or the brief discussion following \cite[Theorem 3.8]{Ske06b}.)
}
Denoting concrete von Neumann modules as pairs $(E,H)$ and representations of $\cB'$ as pairs $(\rho',H)$, the relation
\beqn{
(E,H)
~\longleftrightarrow~
(\rho',H)
}\eeqn
establishes a one-to-one functor between the category of concrete von Neumann \nbd{\cB}modules and the category of normal nondegenerate representations of $\cB'$; see \cite[Theorem 3.10]{Ske06b}. (The morphisms are $\sB^a(E_1,E_2)=\sB^{bil}(H_1,H_2)=C_{\cB'}(\sB(H_1,H_2))$. In particular, $\sB^a(E)=\rho'(\cB')'\subset\sB(H)$.)

If, for another von Neumann algebra $\cA\subset\sB(K)$, a (concrete) von Neumann \nbd{\cB}module $E$ is an \nbd{\cA}\nbd{\cB}correspondence (see the conventions) such that the canonical homomorphism is normal, then we say $E$ is a (\hl{concrete}) \hl{von Neumann correspondence}. For a von Neumann \nbd{\cA}\nbd{\cB}correspondence $E$ this means $\cA\ni a\mapsto a\odot\id_G\in\sB(E\odot G)$ is a normal representation of $\cA$ on $E\odot G$. (Here, we tacitly identify $a\in\cA$ with its action $a\in\sB^a(E)$, even if $E$ is not faithful; recall that a correspondence is \hl{faithful} if the canonical homomorphism is faithful.) For a concrete von Neumann \nbd{\cA}\nbd{\cB}correspondence $E\subset\sB(G,H)$ this means that the \hl{Stinespring representation} $\rho$ defined by setting $\rho(a)\colon xg\mapsto(ax)g$ of $\cA$ associated with $E$ is normal.%
\footnote{
As normality of homomorphisms between von Neumann or \nbd{W^*}algebras is an intrinsic property depending only on the algebraic structure, normality of the left action is defined as soon as $\sB^a(E)$ is proved to be a von Neumann or \nbd{W^*}algebra. (The way to do this is, though, differently easy.) Another entirely intrinsic iff-criterion is to require that the (CP-)maps $\AB{x,\bullet x}\colon\cA\rightarrow\cB$ are normal for all $x\in E$. See \cite[Lemma 3.3.2]{Ske01} for this and one more equivalence.
}
(That $\rho$ is well-defined follows, of course, from the fact that it is unitarily equivalent to $a\mapsto a\odot\id_G$.)

Given a concrete von Neumann \nbd{\cB}module $(E,H)=(\rho',H)$, a representation $\rho\colon\cA\rightarrow\sB(H)$ maps into $\sB^a(E)=\rho'(\cB')'$ if and only if
\beqn{
\bSB{\rho(\cA),\rho'(\cB')}
~=~
\zero,
}\eeqn
that is, if $\rho'$ and $\rho$ are a \hl{commuting pair} of normal and nondegenerate representations. We see, 
\beqn{
(E,H)
~\longleftrightarrow~
(\rho',\rho,H)
}\eeqn
establishes a one-to-one functor between the category of concrete von Neumann \nbd{\cA}\nbd{\cB}cor\-re\-spond\-ences and the category of commuting pairs of normal nondegenerate representations of $\cB'$ and $\cA$. (The morphisms are $\sB^{a,bil}(E_1,E_2)=C_{\cB'}(\sB(H_1,H_2))\cap C_\cA(\sB(H_1,H_2))$. In particular, $\sB^{a,bil}(E)=\rho'(\cB')'\cap\rho(\cA)'\subset\sB(H)$.)

Now, writing a concrete von Neumann \nbd{\cA}\nbd{\cB}correspondence $E$ as triple $(\rho',\rho,H)$ means fixing the first representation $\rho'$ to be the commutant lifting of $\cB'$ and the second representation $\rho$ to be the Stinespring representation of $\cA$ (both associated with $E$).%
\footnote{
For the choice of terminology, observe that for $E$ being (the strong closure of) Paschke's GNS-correspondence of a normal CP-map $T$ from $\cA$ to $\cB$, the representation $\rho$ \bf{is} the Stinespring representation of $T$ and if $\cB$ is generated by $T(\cA)$, then $\rho'$ is the representation of $\cB'$ occurring in the section with the title ``lifting commutants'' of Arveson's \cite{Arv69}.
}
But (as observed first in Skeide \cite{Ske03c} and Muhly and Solel \cite{MuSo05}), nobody prevents us from exchanging these roles, taking $\rho$ as commutant lifting of $(\cA')'=\cA$ and $\rho'$ as the Stinespring representation of $\cB'$ for the concrete von Neumann \nbd{\cB'}\nbd{\cA'}correspondence
\beqn{
E'
~:=~
C_\cA(\sB(K,H))
}\eeqn
corresponding to the triple $(\rho,\rho',H)$. The relation
\beqn{
(E,H)
~=~
(\rho',\rho,H)
~\longleftrightarrow~
(\rho,\rho',H)
~=~
(E',H)
}\eeqn
establishes a one-to-one functor, called \hl{commutant}, between the category of concrete von Neumann \nbd{\cA}\nbd{\cB}correspondences and the category of concrete von Neumann \nbd{\cB'}\nbd{\cA'}cor\-re\-spond\-ences; see \cite[Theorem 4.4]{Ske06b}.%
\footnote{
The commutant, considered as functor on the category of all concrete von Neumann correspondences (without fixing $\cA$ and $\cB$) is, clearly, self-inverse; it generalizes the commutant of von Neumann algebras by considering a von Neumann algebra a correspondence over itself, the \hl{identity correspondence}, in the obvious way.
}
(While for product systems we have $\cA=\cB$, it is not only clearer in general to keep the algebras separate. In the discussion of left and right dilations we actually need the case with $\cA=\C\subset\sB(\C)$.)

\bob \label{commob}
The von Neumann correspondence $E_t:={_t}\cB$ of the member $\vt_t$ of an \nbd{E_0}semi\-group $\vt$ on $\cB\subset\sB(G)$, seen as concrete von Neumann \nbd{\cB}correspondence $E_t\subset\sB(G,{_t}G)=\sB(G)$, is given by the triple $(\id_{\cB'},\vt_t,G)$. Its commutant, given by the triple $(\vt_t,\id_{\cB'},G)$, is, therefore, the concrete von Neumann \nbd{\cB'}module
\beqn{
E'_t
~=~
\bCB{x'_t\in\sB(G)\colon\vt_t(b)x'_t=x'_tb~(b\in\cB)}
}\eeqn
with left action of $\cB'$ given (via $\id_{\cB'}$) simply by operator multiplication. Moreover, we readily verify that
\beqn{
u'_{t,s}
\colon
x'_t\odot y'_s
~\longmapsto~
x'_ty'_s
}\eeqn
(operator multiplication) turns $E'^\odot=\bfam{E'_t}_{t\in\R_+}$ into a product system, the \hl{intertwiner product system} of $\vt$ as discussed in the introduction (in the second item of the discussion following Theorem \ref{pairedthm}). There, we denoted it by $F'^\odot$ to not confuse it with the commutant of $E^\odot$; after Observation \ref{PScommob}, calling it $E'^\odot$ will be fully justified. ($F'^\odot$, in this section, will rather stand for the strongly continuous commutant of $E^\odot$, following Definition \ref{sccdefi}.)
\eob

The product system in this observation is a special case of the intertwiner product system of an \nbd{E_0}semigroup $\vt$ on $\sB^a(E)$ constructed in Skeide \cite[Section 2]{Ske03c} and recognized (in a very hand-waving way) as being isomorphic to the commutant of \bf{the} product system of that \nbd{E_0}semigroup. (We do not need it, but still give some hints in Footnote \ref{rduniFN}.) As far as the algebraic commutant system is concerned, all we need is Observation \ref{commob} plus the following discussion:

\bulletline[The tensor product under commutant.~]
In Observation \ref{commob}, we have identified the members $E'_t$ of the intertwiner system $E'^\odot$ of an \nbd{E_0}semigroup $\vt$ on $\cB$ as the commutants of the members $E_t$ of \bf{the} product system $E^\odot$ of $\vt$. It is true (and easily verified) that the product $u'_{t,s}$ as indicated does define a product system structure on the family formed by the $E'_t$. But is this product system structure indeed the one inherited from being the commutant of $E^\odot$?

So far we did believe the (positive) answer to this question from \cite[Section 2]{Ske03c}. But to be specific about several things -- more so than the hand-waiving discussion in \cite{Ske03c} -- we need to understand the tensor product under commutant more specifically, too. 

Let us recall that for von Neumann algebras $\cA\subset\sB(K)$, $\cB\subset\sB(G)$, and $\cC\subset\sB(L)$ and concrete von Neumann correspondences $E=(\rho',\rho,H)$ and $F=(\pi',\pi,\wt{H})$ from $\cA$ to $\cB$ and from $\cB$ to $\cC$, respectively, their \hl{tensor product} is $E\sodots F=(\id_E\odot\pi',\rho\odot\id_{\wt{H}},E\odot\wt{H})$. The commutant
\beqn{
(E\sodots F)'
~=~
(\rho\odot\id_{\wt{H}},\id_E\odot\pi',E\odot\wt{H})
}\eeqn
is canonically isomorphic to $F'\sodots E'=(\id_{F'}\odot\rho,\pi'\odot\id_H,F'\odot H)$ via
\beqn{
x\odot y'g
~\longmapsto~
y'\odot xg,
}\eeqn
taking also into account that $\cls F'G=\wt{H}$ and $\cls EG=H$. (This is easier to understand looking at the picture $E\odot F'\odot G\cong F'\odot E\odot G$ via $x\odot y'\odot g\mapsto y'\odot x\odot g$, clearly isometric with total range, and how the several representations are intertwined.)%
\footnote{
The tensor product of von Neumann correspondences and \nbd{W^*}correspondences both start with $E\odot F$ (tensor product of \nbd{C^*}correspondences) and, then, close it suitably. While for von Neumann correspondences this closure it plain (simply take the strong closure in $\sB(L,(E\odot F)\odot L)$), for \nbd{W^*}correspondences this is much more cumbersome. Already indicating the elements of the space, namely of $\sB^r(E\odot F,\cC)$, is not so easy; but it gets really nasty when we wish to concretely calculate their inner products. (Following Paschke \cite{Pas73}, one has to approximate at least one of them by elements of $E\odot F$. Already, Rieffel \cite{Rie74a} noticed the possibility to obtain the space as strong closure in a suitable operator space where the computation of inner products is easy with the help of the adjoint as operator multiplication.)

Note that the tensor product of two concrete von Neumann correspondences is ``the same'' as their tensor product when considered as von Neumann correspondences by the canonical isomorphism of $E\odot(F\odot L)$ and $(E\odot F)\odot L$.

Note, too, that the tensor product of concrete von Neumann correspondences includes the \it{fusion product} of \it{Connes correspondences} as the special case when all von Neumann algebras come in standard representation; see \cite[Section 3]{Ske22b}
}

Recall that a morphism $a\in\sB^{a,bil}(E\sodots F,D)\subset\sB(E\odot\wt{H},\wh{K})$ into another concrete von Neumann \nbd{\cA}\nbd{\cC}correspondence $D=(\sigma',\sigma,\wh{K})$ gives rise to a morphism $(E\sodots F)'\rightarrow D'$ represented by the same $a\in\sB^{a,bil}((E\sodots F)',D')=\sB^{a,bil}(E\sodots F,D)\subset\sB(E\odot\wt{H},K)$. Under the canonical isomorphism $E\odot\wt{H}\cong F'\odot H$, this induces a morphism $a'\in\sB^{a,bil}(F'\sodots E',D')$.

\bob \label{PScommob}
In the notations of Observation \ref{commob}: $u'_{t,s}\colon E'_t\sodots E'_s\rightarrow E'_{t+s}$ is, indeed, the isomorphism $(u_{s,t})'$ induced from $u_{s,t}\colon E_s\sodots E_t\rightarrow E_{s+t}$ from taking commutants. Namely:
\beqn{
E'_t\sodots E'_s
~=~
(\id_{E'_t}\odot\vt_s,\id_{\cB'}\odot\id_{_sG},E'_t\odot {_s}G),
}\eeqn
while
\beqn{
(E_s\sodots E_t)'
~=~
(\vt_s\odot\id_{_tG},\id_{E_s}\odot\id_{\cB'},E_s\odot{_t}G).
}\eeqn
The two are canonically isomorphic by identifying $E'_t\odot {_s}G$ with $E_s\odot{_t}G$ via
\beqn{
x'_t\odot(y_sg)
~\longmapsto~
y_s\odot(x'_tg)
}\eeqn
(carrying along also the right representations). Plugging in on the right-hand side a typical element $g_t=x_tg$ of $_tG=G$, we see that $u_{s,t}(y_s\odot g_t)=\vt_t(y_s)g_t\in{_{s+t}}G=G$. Especially,
\beqn{
u_{s,t}(y_s\odot(x'_tg))
~=~
\vt_t(y_s)(x'_tg)
~=~
x'_ty_sg.
}\eeqn
Plugging in on the left-hand side a typical element $g'_s=y'_sg$ of $_sG=G$, we see that $u'_{t,s}(x'_t\odot g'_s)=x'_tg'_s\in{_{s+t}}G=G$. Especially,
\beqn{
u'_{t,s}(x'_t\odot(y_sg))
~=~
x'_t(y_sg)
~=~
x'_ty_sg.
}\eeqn
We see, the two coincide, that is, the $(u_{s,t})'$ induced on the commutants from $u_{s,t}$ is $u'_{t,s}$ as we defined it.
\eob

\bulletline[Right dilations and representation.~]
\it{Right dilations} occurred, just as a name, in Section \ref{pairSEC}, and finding a \it{strongly continuous right dilation} is an essential step in Outline \ref{outob} for the construction of the strongly continuous commutant. We arrive at the notion of right dilation most easily, taking into account the behaviour of the commutant we just discussed, as commutant of a left dilation.

Suppose $v_t$ is a left dilation of $E^\odot$ to the concrete von Neumann \nbd{\cB}module $(E,H)=(\rho',H)$ and consider $E$ as concrete von Neumann \nbd{\C}\nbd{\cB}correspondence $(\rho',\id_H\id_\C,H)$ for the von Neumann algebra $\C=\sB(\C)$, with the Stinespring representation $\id_H\id_\C\colon z\mapsto\id_Hz$ of $\C$ on $H$. Then the commutant of this concrete von Neumann \nbd{\C}\nbd{\cB}cor\-re\-spond\-ence $(\rho',\id_H\id_\C,H)$ is the concrete von Neumann \nbd{\cB'}\nbd{(\C'\!\!\!=\!\!\C)}correspondence $(\id_H\id_\C,\rho',H)$, that is, the Hilbert space $H=\sB(\C,H)$ with the left action of $\cB'$ via $\rho'$. The family of \nbd{\C}\nbd{\cB}linear unitaries $v_t\colon E\sodots E_t\rightarrow E$ gives rise to a family of \nbd{\cB'}\nbd{\C}linear unitaries $w'_t\colon E'_t\odot H\rightarrow H$ (since $H$ is a Hilbert space, no strong closure is necessary) which iterate associatively with the product system structure of $E'^\odot$. Moreover, $H$ is faithful%
\footnote{
It is easy to see that a general concrete von Neumann correspondence is faithful if and only if its commutant is strongly full; and the $E$ of a left dilation to $E$ is, by definition, strongly full.
}
and we recover $\vt$ as $\vt_t=w'_t(\id'_t\odot\bullet)w'^*_t$ (where, recall!, $\sB^{bil}(H)=\rho'(\cB')'=\sB^a(E)$ is the algebra on which the right-hand side naturally acts).

This motivates the definition part (adding also strong continuity) of the following.

\bitemp[Definition and Theorem {\cite[Appendix B.2]{Ske16}}.~] \label{rdthm}
A \hl{right dilation} of a (necessarily faithful) product system $E^\odot$ of (concrete) von Neumann \nbd{\cB}correspondences to a Hilbert space $H$ with a faithful normal nondegenerate representation $\rho$ of $\cB$ (that is, to a faithful von Neumann \nbd{\cB}\nbd{\C}correspondence $H$) is a family of bilinear unitaries $w_t\colon E_t\odot H\rightarrow H$ that iterate associatively with the product of $E^\odot$. We use the product notation $x_th:=w_t(x_t\odot h)$.

A right dilation $w_t$ of a strongly continuous product system $E^\odot$ to $H$ is strongly continuous if the function $t\mapsto x_th\in H$ is continuous for all $\bfam{x_t}_{t\in\R_+}\in CS_i^s(E^\odot)$.

By \cite[Theorem B.31]{Ske16}, every faithful strongly continuous product system admits a strongly continuous right dilation.
\eitemp

As in the discussion leading to that definition, for every right dilation $w_t$, the maps $\vt^w_t:=w_t(\id_t\odot\bullet)w_t^*$ define an \nbd{E_0}semigroup $\vt^w$ (strongly continuous, if the right dilation is strongly continuous) on $\sB^{bil}(H)=\rho(\cB)'=\sB^a(E')$ (where we consider the concrete von Neumann \nbd{\cB'}module $(\rho,H)=(E',H)$). Moreover the intertwiner system of $\vt^w$ (as defined in \cite[Section 2]{Ske03c}) is (as discussed there) isomorphic to the commutant of $E^\odot$.%
\footnote{ \label{rduniFN}
The discussion preceding Definition \ref{rdthm}, actually, furnishes a one-to-one correspondence between left and right dilations of a fixed product system and its commutant, respectively, and the left and right dilation in such a pair determine the same \nbd{E_0}semigroup. Since by the last item of Observation \ref{TPSob} the \nbd{E_0}semigroup constructed from a left dilation of a product system determines that product system up to isomorphism to be \bf{the} product system of the \nbd{E_0}semigroup, the same statement is true for the corresponding right dilation and the intertwiner system. This simple observation is the reason that allowed us in our earlier papers to be more unspecific about the concrete choice inside the isomorphism class: The intertwiner system of an \nbd{E_0}semigroup comes along with a right dilation, thus assuring that the intertwiner system is isomorphic to the commutant of \bf{the} product system. Here, however, to identify the candidates for the strongly continuous sections of intertwiners we have to be specific.
}

\vspace{1ex}
This concludes the the discussion of all missing pieces in Outline \ref{outob}:

\vspace{-1ex}\bdefi \label{sccdefi}
The \hl{strongly continuous commutant} of a faithful strongly continuous product system is the strongly continuous product system constructed following Outline \ref{outob}.
\edefi

By (the proof of) \cite[Theorem B.34]{Ske16}, the strongly continuous commutant of $E^\odot$ does not depend on the choice of the strongly continuous right dilation and, if $E^\odot$ is also strongly full, the strongly continuous double commutant gives back $E^\odot$.

\lf
There is a concept to look at right dilations, called \it{representation}, which frequently comes in usefully for notational reasons.

\bdefi \label{PSrepdefi}
A \hl{representation} of a product system $E^\odot$ (of von Neumann correspondences) on a Hilbert space $H$ is family of maps $\eta_t\colon E_t\rightarrow\sB(H)$ such that
\baln{
\eta_t(x_t)\eta_s(y_s)
&
~=~
\eta_{t+s}(x_ty_t),
&
\eta_t(x_t)^*\eta_t(y_t)
&~=~
\eta_0(\AB{x_t,y_t)}.
}\ealn
It follows that $\pi:=\eta_0$ is a (normal) representation of $\cB$ and that each $\eta_t$ is a bimodule map with respect to $\pi$. A representation is \hl{normal} if $\pi$ is normal (so that each $\eta_t$ is \nbd{\sigma}weak). A representation is \hl{faithful} if $\pi$ (and, therefore, each $\eta_t$) is faithful. A representation is \hl{nondegenerate} if each $\eta_t$ acts nondegenerately on $H$.
\edefi

\it{Representation} goes under the name \it{isometric covariant representation} in Muhly and Solel \cite{MuSo02}. The concept of \it{faithful normal nondegenerate representation} $\eta_t$ is equivalent to the concept of \it{right dilation} via the one-to-one correspondence determined by
\beqn{
\eta_t(x_t)
~=~
w_t(x_t\odot h)
}\eeqn
(\cite[Proposition B.19]{Ske16}). Nondegeneracy of $\eta_t$ corresponds to surjectivity of $w_t$.

\bulletline[Strongly continuous intertwiner sections.~]
Returning to our problem, we now fix $\vt$ to be a strongly continuous and faithful \nbd{E_0}semigroup on the von Neumann algebra $\cB\subset\sB(G)$. \bf{The} product system of $\vt$, $E^\odot$, as pointed out in Observation \ref{TPSob}, is the one constructed from the unit vector $\U\in\cB$. The embeddings $i_t$ are just the identification maps $_t\cB=\cB$ as right modules. The strongly continuous sections are, therefore, just the strongly continuous functions $\R_+\rightarrow\cB$. 

We use the notations from Observations \ref{commob} and \ref{PScommob}, in which we identified very carefully the intertwiner system of $\vt$, $E'^\odot$, as the (algebraic) commutant system of $E^\odot$. We shall now construct the strongly continuous commutant $F'^\odot$ of $E^\odot$ as outlined in \ref{outob}, and then give an explicit (algebraic; there is no strongly continuous structure on $E'^\odot$) isomorphism of $E'^\odot$ with $F'^\odot$. Only then, we show that the images of the strongly continuous sections of $F'^\odot$ in $E'^\odot$ act strongly continuously on $G$ (Lemma \ref{iclem}) leading as a corollary to Theorem \ref{ecthm} stating that the intertwiner system has enough strongly continuous sections.

As a strongly full and faithful strongly continuous product system, $E^\odot$ has a strongly full and faithful strongly continuous commutant, which, by Definition \ref{sccdefi}, can be obtained as outlined in \ref{outob} in the following way:
\begin{itemize}
\item
Take any strongly continuous right dilation $w_t\colon E_t\odot H\rightarrow H$ of $E^\odot$ to a Hilbert space $H$ with a (by definition of right dilation) faithful normal nondegenerate representation $\rho$ of $\cB$. (Recall that by \cite[Theorem B.31]{Ske16} such a dilation always exists.)

For notational convenience, we also define the (nondegenerate faithful) representation $\eta_t(x_t)\colon h\mapsto w_t(x_t\odot h)$ of $E^\odot$ on $H$. (Recall that $\eta_0=\rho$.) Then the condition that $w_t$ be strongly continuous, means exactly that $t\mapsto\eta_t(x_t)\in\sB(H)$ is strongly continuous for every strongly continuous section $\bfam{x_t}_{t\in\R_+}$ of $E^\odot$.

\item
Let $E':=C_\cB(\sB(G,H))$ denote the concrete von Neumann \nbd{\cB'}module given by the pair $(\rho,H)$. Recall that $\sB^a(E')=\rho(\cB)'\subset\sB(H)$.

\item
Define the \nbd{E_0}semigroup $\vt^w$ on $\rho(\cB)'=\sB^a(E')$ as $\vt^w_t=w_t(\id_t\odot\bullet)w_t^*$. Since $w_t$ is strongly continuous, so is $\vt^w$.

\item
Assume that $E'$ has a unit vector $\xi'$. (By \cite[Lemma 4.2]{Ske09}, we may achieve this by amplifying $w_t$ to $H^\en$ for a suitable cardinality $\en$.)

Recall that $b'\mapsto\xi'b'\xi'^*$ defines a (usually, degenerate) normal faithful representation $\xi'\bullet\xi'^*$ of $\cB'$ on $H$ that acts as $\xi'b'\xi'^*\colon\xi'g\mapsto\xi'b'g$ on the subspace $\xi'G$ of $H$ (isomorphic to $G$ via the isometry $\xi'$) and as $0$ on its complement. We may say $G$ sits in $H$, and $\cB'$ acts, there, as it should.

\item
The strongly continuous commutant is, then, the product system $F'^\odot$ of \nbd{\cB'}cor\-re\-spond\-ences $F'_t$ constructed from the strongly continuous  \nbd{E_0}semigroup $\vt^w$ on $\sB^a(E')$ with the help of the unit vector $\xi'$. In other words:
\begin{itemize}
\item
$F'_t:=\vt^w_t(\xi'\xi'^*)E'\subset\sB(G,H_t)\subset\sB(G,H)$ (where we defined the subspace $H_t:=\vt^w_t(\xi'\xi'^*)H$ of $H$) as right submodule of $E'$.

\item
The left action of $\cB'$ on $F'_t$ is given by $\vt^w_t(\xi'\bullet\xi'^*)$.

\item
This means $F'_t$, as concrete von Neumann correspondence, is given by the triple $(\rho\upharpoonright H_t\,,\,\vt^w_t(\xi'\bullet\xi'^*)\upharpoonright H_t\,,\,H_t)$.

\item
The product of $F'^\odot$ is given by $x'_s\odot y'_t\mapsto\vt^w_t(x'_s\xi'^*)y'_t$.

\item
The strongly continuous sections of $F'^\odot$ are those sections such that the map $t\mapsto x'_t\in F'_t\subset E'\subset\sB(G,H)$ is strongly continuous.
\end{itemize}
Recall that for each $s$ and each $x'_s=\vt^w_s(\xi'\xi'^*)x'\in F'_s$ the section $\bfam{\vt^w_t(\xi'\xi'^*)x'}_{t\in\R_+}$ is strongly continuous and meets $x'_s$ at $t=s$.
\end{itemize}

\noindent
We now show that $F'^\odot$ is (algebraically) isomorphic to $E'^\odot$ in a specific way.

Recall that $H=\cls E'G$ and that $H=\cls\eta_t(E_t)H$, by nondegeneracy of $\eta_t$, so that
\beqn{
H
~=~
\cls\eta_t(E_t)E'G.
}\eeqn
Since $\vt^w_t(a')\eta_t(x_t)h=w_t(x_t\odot a'h)=\eta_t(x_t)a'h$, we find
\beq{ \label{etavtwint}
\vt^w_t(a')\eta_t(x_t)
~=~
\eta_t(x_t)a',
}\eeq
hence,
\beq{ \label{Htetattot}
H_t
~=~
\vt^w_t(\xi'\xi'^*)\cls\eta_t(E_t)E'G
~=~
\cls\eta_t(E_t)\xi'\xi'^*E'G
~=~
\cls\eta_t(E_t)\xi'G.
}\eeq
By
\bmun{
\AB{\eta_t(x_t)\xi'g,\eta_t(y_t)\xi'\hat{g}}
~=~
\AB{\xi'g,\eta_t(x_t)^*\eta_t(y_t)\xi'\hat{g}}
\\
~=~
\AB{\xi'g,\rho(\AB{x_y,y_t})\xi'\hat{g}}
~=~
\AB{\xi'g,\xi'\AB{x_y,y_t}\hat{g}}
~=~
\AB{g,\AB{x_y,y_t}\hat{g}},
}\emun
we see that the map $\Upsilon_t$ defined by
\beqn{
\Upsilon_t
\colon
G\ni x_tg
~\longmapsto~
\eta_t(x_t)\xi'g\in H_t
}\eeqn
is a unitary $G\rightarrow H_t$. Repeating the computation in \eqref{Htetattot} with $\vt^w_t(\xi'b'\xi'^*)$, we see from
\beqn{
\vt^w_t(\xi'b'\xi'^*)\bfam{\eta_t(x_t)\xi'g}
~=~
\eta_t(x_t)\xi'b'g
}\eeqn
that $\Upsilon_t$ intertwines the actions of $\cB'$. Since
\beqn{
\rho(b)\bfam{\eta_t(x_t)\xi'g}
~=~
\eta_t(b.x_t)\xi'g,
}\eeqn
$\Upsilon_t$ also intertwines the actions of $\cB$. Therefore, since $\Upsilon_t$ sends the triple $(\vt_t,\id_{\cB'},G)$ to the triple $(\rho\upharpoonright H_t\,,\,\vt^w_t(\xi'\bullet\xi'^*)\upharpoonright H_t\,,\,H_t)$ (intertwining the occurring pairs of representations), the map $\upsilon_t\colon E'_t\ni x'_t\mapsto \Upsilon_tx'_t\in F'_t$ is an isomorphism $\sB(G)\supset E'_t\rightarrow F'_t\subset\sB(G,H_t)$.

\bprop \label{isoprop}
The $\upsilon_t$ form an isomorphism of product systems $E'^\odot\rightarrow F'^\odot$.
\eprop

\proof
The product of elements $x'_s,\hat{x}'_t$ in $E'^\odot$ is simply $x'_s\hat{x}'_t$. The product of the elements $\upsilon_sx'_s,\upsilon_t\hat{x}'_t$ in $F'^\odot$ is $\vt^w_t((\upsilon_sx'_s)\xi'^*)(\upsilon_t\hat{x}'_t)$. We have to show that the latter coincides with $\upsilon_{s+t}(x'_s\hat{x}'_t)$.

We will check this weakly (recalling that the product of elements $x_t,\hat{x}_s$ in $E^\odot$ is $\vt_s(x_t)\hat{x_s}$ and that the $\eta_t$ form a representation of $E^\odot$) with elements
\beqn{
Z
~:=~
\Upsilon_{s+t}\vt_s(x_t)\hat{x_s}g
~=~
\eta_t(x_t)\eta_s(\hat{x}_s)\xi'g,
}\eeqn
which are total in $H_{s+t}$. We start with $\upsilon_{s+t}(x'_s\hat{x}'_t)$ and find
\bmu{ \label{ust}
\AB{Z,\upsilon_{s+t}(x'_s\hat{x}'_t)\hat{g}}
~=~
\AB{\Upsilon_{s+t}\vt_s(x_t)\hat{x_s}g,\Upsilon_{s+t}x'_s\hat{x}'_t\hat{g}}
~=~
\AB{\hat{x_s}g,\vt_s(x_t^*)x'_s\hat{x}'_t\hat{g}}
~=~
\AB{{x'_s}^*\hat{x_s}g,x_t^*\hat{x}'_t\hat{g}},
}\emu
where we made use of the fact that $x'_s$ intertwines $\vt_t(b)$ (for any element $b\in\cB$, hence, also for the element $x_t^*\in\cB$). For $\vt^w_t((\upsilon_sx'_s)\xi'^*)(\upsilon_t\hat{x}'_t)$ we find
\bmu{ \label{usut}
\AB{Z,\vt^w_t((\upsilon_sx'_s)\xi'^*)(\upsilon_t\hat{x}'_t)\hat{g}}
~=~
\AB{\eta_t(x_t)\eta_s(\hat{x}_s)\xi'g,\vt^w_t((\upsilon_sx'_s)\xi'^*)(\upsilon_t\hat{x}'_t)\hat{g}}
\\
~=~
\AB{\vt^w_t(\xi'(\upsilon_sx'_s)^*)\eta_t(x_t)\eta_s(\hat{x}_s)\xi'g,\Upsilon_t\hat{x}'_t\hat{g}}
~=~
\AB{\Upsilon_t^*\eta_t(x_t)\xi'(\upsilon_sx'_s)^*\eta_s(\hat{x}_s)\xi'g,\hat{x}'_t\hat{g}},
}\emu
where we made use of \eqref{etavtwint}. Since $\Upsilon_tx_t\colon g\mapsto\eta(x_t)\xi'g$, we see that $\Upsilon_t^*\eta_t(x_t)\xi'=x_t$. Continuing from \eqref{usut}, we obtain
\bmun{
\AB{Z,\vt^w_t((\upsilon_sx'_s)\xi'^*)(\upsilon_t\hat{x}'_t)\hat{g}}
~=~
\AB{(\upsilon_sx'_s)^*\eta_s(\hat{x}_s)\xi'g,x_t^*\hat{x}'_t\hat{g}},
\\
~=~
\AB{{x'_s}^*\Upsilon_s^*\eta_s(\hat{x}_s)\xi'g,x_t^*\hat{x}'_t\hat{g}}
~=~
\AB{{x'_s}^*\hat{x_s}g,x_t^*\hat{x}'_t\hat{g}}.
}\emun
Confronting this with \eqref{ust}, se get $\upsilon_{s+t}(x'_s\hat{x}'_t)=\vt^w_t((\upsilon_sx'_s)\xi'^*)(\upsilon_t\hat{x}'_t)$.\qed

\blem \label{iclem}
The inverse isomorphism $\upsilon_t^*$ sends strongly continuous sections $\bfam{x'_t}_{t\in\R_+}$ of the strongly continuous product system $F'^\odot$ to sections $\bfam{\upsilon_t^*x'_t}_{t\in\R_+}$ of the intertwiner product system $E'^\odot$ which are strongly continuous as maps $t\mapsto \upsilon_t^*x'_t\in\sB(G)$.
\elem

\proof
Let $\bfam{x'_t}_{t\in\R_+}$ be in $F'^\odot$ as stated.

We first show weak continuity of $y'_t:=\upsilon_t^*x'_t$. We have to show that
\beqn{
t
~\longmapsto~
\AB{g,\Upsilon_t^*x'_t\hat{g}}
~=~
\AB{\Upsilon_tg,x'_t\hat{g}}
}\eeqn
is continuous for all $g,\hat{g}\in G$. Since the section of $F'^\odot$ is strongly continuous, the function $t\mapsto x'_t\hat{g}\in H_t\subset H$ is continuous for any $\hat{g}\in G$.  Putting $x_t:=\U_t\in E_t={_t}\cB=\cB$, we have $x_tg=g$ for all $t$, so that
\beqn{
\Upsilon_tg
~=~
\Upsilon_tx_tg
~=~
\eta_t(x_t)\xi'g.
}\eeqn
Since the section $\bfam{x_t}_{t\in\R_+}$ is strongly continuous and since $\eta_t$ is strongly continuous, also the function $t\mapsto\Upsilon_tg\in H_t\subset H$ is continuous for any $g\in G$. In conclusion, $t\mapsto y'_t\in E'_t\subset\sB(G)$ is weakly continuous.

The fact that this function is actually strongly continuous, follows from the well-known observation that $y'_tg\to y'_sg$ in norm if (and only if) $y'_tg\to y'_sg$ weakly and $\snorm{y'_tg}\to\snorm{y'_sg}$, and from the fact that $\snorm{x'_tg}\to\snorm{x'_sg}$ for the corresponding section $x'_t$ in $F'^\odot$ (which \bf{is} strongly continuous) and that $\snorm{y'_tg}=\snorm{x'_tg}$.\qed

\lf
Theorem \ref{ecthm} is now a simple corollary of Lemma \ref{iclem}:

\proof[Proof of Theorem \ref{ecthm}.~]
Since the strongly continuous sections of $E'^\odot$ constructed in Lemma \ref{iclem} are exactly the images of the strongly continuous sections of $F'^\odot$ under the (inverse of the) isomorphism in Proposition \ref{isoprop}, and since $F'^\odot$, as strongly continuous product system, admits enough strongly continuous sections, $E'^\odot$ admits enough strongly continuous sections, too.\qed

\section{Proof of Theorem \ref{cpairedthm}} \label{cpairSEC}

Recall that we are back now to the notations from Section \ref{pairSEC} (as described also in the introduction) and that the results from Section \ref{E0psSEC} (especially, the isomorphism from Proposition \ref{isoprop}) are applied to $\vt'$ rather than to $\vt$.

As emphasized in the discussion leading to Theorem \ref{ipairedthm}, in Section \ref{pairSEC} we have not just proved Theorem \ref{pairedthm}, but we, actually,  have established a specific one-to-one correspondence between unitary semigroups (modulo a multiplier) $U_t$ implementing a paring, and isomorphisms (modulo the same multiplier) $u_t$ between the product system of one \nbd{E_0}semigroup and the commutant of the product system of the other \nbd{E_0}semigroup. With a possible multiplier we get Theorem \ref{pairedthm}; without a multiplier we get Theorem \ref{ipairedthm}.

By Wigner's theorem, strongly continuous pairings \bf{are} inner (so, we are in the situation of Theorem \ref{ipairedthm}) and the implementing unitary (semi)group \bf{can} be chosen strongly continuous. Theorem \ref{ipairedthm} leaves us with the question whether, for strongly continuous $\vt$ and $\vt'$\!, strong continuity of the unitary group implementing an inner paring is reflected exactly by strong continuity of the corresponding (algebraic) isomorphism of strongly continuous product systems.

\proof[Proof of Theorem \ref{cpairedthm} ''$\Rightarrow$''.~]
Suppose the unitary group $U_t$ implementing a pairing $\alpha$ is strongly continuous. Recall that a section $\bfam{x_t}_{t\in\R_+}$ in $E^\odot=\bfam{_t\cB}_{t\in\R_+}$ is strongly continuous, if $t\mapsto x_t\in E_t\subset\sB(G)$ is strongly continuous (as map into $\sB(G)$). Therefore, the map
\beqn{
t
~\longmapsto~
U_tx_t
~\in~
F_t
~\subset~
\sB(G)
}\eeqn
is strongly continuous (as map into $\sB(G)$), too.

Recall that, in order to obtain the strongly continuous commutant of $E^\odot$, we may take just any strongly continuous right dilation $w_t\colon E_t\odot H\rightarrow H$ of $E^\odot$ to a von Neumann \nbd{\cB}\nbd{\C}cor\-re\-spond\-ence (that is, a Hilbert space with normal nondegenerate left action of $\cB$) such that the von Neumann \nbd{\cB'}module $E':=C_\cB(\sB(G,H)$ admits a unit vector $\xi'$, and that the strongly continuous commutant is, then, given by the strongly continuous product system of the strongly continuous \nbd{E_0}semigroup $\vt^w_t:=w_t(\id_t\odot\bullet)w_t^*$ on $\sB^{bil}(H)=\sB^a(E')$.

Our choice for this right dilation is just $H=G$ (with the identity representation of $\cB$) and the right dilation $w_t(x_t\odot g):=U_tx_tg$.%
\footnote{
It is easy to check directly that this is a right dilation. A more structural way to argue, is the following (easy to check) insight: If we have an isomorphism $u_t$ from $E^\odot$ to $F^\odot$ and a right dilation $\hat{w}_t$ of $F^\odot$ to $H$, then $w_t:=\hat{w}_t(u_t\odot\id_H)$ is a right dilation of $E^\odot$. (In our case, we take for $\hat{w}_t$ the identity right dilation of the intertwiner system $F^\odot$.)
}
By the first paragraph, this right dilation is strongly continuous. Moreover, $E'=\cB'$ contains a unit vector $\xi'=\U'$ and, applying $\vt^w_t(b')$ to elements $U_tx_tg=w_t(x_t\odot g)$, we get
\beqn{
\vt^w_t(b')U_tx_tg
~=~
w_t(\id_t\odot b')(x_t\odot g)
~=~
w_t(x_t\odot b'g)
~=~
U_tx_tb'g
~=~
\vt'_t(b')U_tx_tg,
}\eeqn
so that $\vt^w$ is just $\vt'$ on $\sB^a(\cB')=\cB'$.

We see, the strongly continuous commutant of the strongly continuous product system of $\vt$ is the strongly continuous product system of $\vt'$.\qed

\lf
Since the most difficult part of the backwards direction (Lemma \ref{iclem}) has been dealt with in Section \ref{E0psSEC}, here we can be rather short.

\proof[Proof of Theorem \ref{cpairedthm} ''$\Leftarrow$''.~]
Assume we have a strongly continuous isomorphism $\hat{u}_t$ from the product system $E^\odot$ of $\vt$ to the strongly continuous commutant of the strongly continuous product system of $\vt'$. By Proposition \ref{isoprop} we have an (algebraic) isomorphism $\hat{\hat{u}}_t$ from that strongly continuous commutant to the intertwiner system $F^\odot$ of $\vt'$, and by Lemma \ref{iclem} $\hat{\hat{u}}_t$ sends strongly continuous sections of the former to sections in the intertwiner system that act strongly continuously. Altogether, $u_t:=\hat{\hat{u}}_t\hat{u}_t$ is an (algebraic) isomorphism that sends strongly continuous sections of $E^\odot$ to strongly continuous intertwiner sections of $F^\odot$.

Define the unitary (semi)group $U_t$ as in Equation \ref{Utdef}, which implements an inner pairing $\alpha$.  Since the $U_t$ are bounded uniformly, for showing strong continuity of $U_t$ in $t=s$ ($s$ fixed) it is enough to show strong continuity on the total subset $F_sG$ of $G$. Every $y_s\in F_s$ can be written $u_sx^0_s$ for a (unique) $x^0_s\in E_s$. There exists a (bounded by $\norm{y_s}$, if we wish) strongly continuous section $\bfam{x_t}_{t\in\R_+}$ of $E^\odot$ such that $x_s=x_s^0$.
\bmun{
U_t(x_s^0g)-U_s(x_s^0g)
~=~
\bfam{U_t(x_sg)-U_t(x_tg)}+\bfam{U_t(x_tg)-U_s(x_sg)}
\\
~=~
U_t((x_sg)-(x_tg))+\bfam{(u_tx_t)g-(u_sx_s)g}.
}\emun
$x_tg$ converges to $x_sg$, and since $\norm{U_t}=1$, so does the first summand. By what we said in the first paragraph, the $u_tx_t$ act strongly continuously on $G$, so, $(u_tx_t)g$ converges to $(u_sx_s)g$. So, $U_t$ converges to $U_s$ strongly on the total subset $F_sG$, hence, everywhere.\qed

\section{Towards a symmetric theory?} \label{scsSEC}

(See also Remark \ref{rem2nd}.) By Theorem \ref{ecthm}, the intertwiner system $E'^\odot$ of a strongly continuous faithful \nbd{E_0}semigroup $\vt$ on $\cB\subset\sB(G)$ has enough sections that act strongly continuously on $G$. But that does not mean that the $E'_t$ sitting in $\sB(G)$, form a strongly continuous product system for the canonical embedding $E'_t\rightarrow\sB(G)$. In fact, the commutant lifting $\vt_t$ of $E'_t$ does depend on $t$, and since the $\vt_t$ are unital (so that all $E'_t$ act nondegenerately on $G$), there cannot be submodules of a single concrete von Neumann \nbd{\cB'}module $(\rho,G)$ (for a fixed commutant lifting $\rho$ of the commutant of $\cB'$), as required in the definition of strongly continuous product system.

However, except for the fact that the canonical embeddings $i'_t\colon E'_t\rightarrow\sB(G)$ do not go into a fixed von Neumann \nbd{\cB'}module, all other requirements of the definition of strongly continuous product system are fulfilled. It would be nice to have a more flexible definition of strongly continuous product system that covers both \bf{the} strongly continuous product system (with its continuous structure coming from a unit vector) and the intertwiner system (with the continuous structure coming from being subspaces of the given $\sB(G)$ as just indicated).

\brem
Note, however, how easy it was to obtain the strongly continuous structure for \bf{the} product system of a strongly continuous \nbd{E_0}semigroup, and how hard it was to get it for the intertwiner system not before running through the whole theory exposed in \cite[Appendix B]{Ske16}. And we had to add to it the missing piece Theorem \ref{ecthm}, a consequence of \cite{Ske16}, before even dreaming of being able to take the strongly constinuous structure from $E'_t$ sitting in $\sB(G)$.
\erem

What we will have in such a setting is a fixed Hilbert space $H$ and a family of commuting pairs of normal representations $\rho_t$ and $\rho'_t$ of $\cB$ and $\cB'$, respectively, on $H$. They will, in general, not be nondegenerate but such that $\rho(\U)=:p_t:=\rho'(\U')$ so that with $H_t:=p_tH$ we have the concrete von Neumann correspondences $E_t:=(\rho'_t,\rho_t,H_t)$ and their commutants $E'_t:=(\rho_t,\rho'_t,H_t)$. We will require that there is a unitary $\Xi\colon G\rightarrow H_0$ such that $\rho^{(')}_0(b^{(')})\Xi=\Xi b^{(')}$, identifying $E_0$ with $\cB$ and $E'_0$ with $\cB'$  via $\Xi$. We know, if there are isomorphisms $u_{s,t}$ turning the $E_t$ into a product system $E^\odot$, then the $u'_{t,s}:=(u_{s,t})'$ turn the $E'_t$ into a product system $E'^\odot$ -- and \it{vice versa}. So far, this is a symmetric setting of embedded algebraic product systems of concrete von Neumann correspondences and their commutants. (In fact, we have the canonical embeddings $i^{(')}_t\colon\sB(G,H_t)\supset E^{(')}_t\rightarrow\sB(G,H)$.)

We will certainly wish that the sets of strongly continuous sections of these product systems are invariant under the bimodule actions of $\cB$ and of $\cB'$, respectively. This can be achieved by requiring that the maps $t\mapsto\rho_t$ and $t\mapsto\rho'_t$ are strongly continuous. And we
wish that we have enough strongly continuous sections of $E^\odot$ and that the product is strongly continuous in the sense that
\beqn{
t
~\longmapsto~
i_{s+t}(x_sy_t)g
}\eeqn
is continuous for all strongly continuous sections (meaning $E^\odot$ is \hl{strongly continuous} in this new sense), and we will wish the analogue conditions for $E'^\odot$ (meaning $E'^\odot$ is \hl{strongly continuous} in this new sense).

It is unclear to what extent these conditions are redundant. Certainly, we would wish to know several things, for now always assuming the condition of strong continuity in $t$ for the representations $\rho_t$ and $\rho'_t$:

\begin{itemize}
\item
Does $E^\odot$ having enough strongly continuous sections, imply already the same statement for $E'^\odot$? And if yes, will also the product of $E'^\odot$ be strongly continuous?

\item
Is $E^\odot$ being strongly full (that is, all $\rho'_t$ are faithful) and strongly continuous in the new sense enough for showing existence of a strongly continuous left dilation (concluding the relation between \nbd{E_0}semigroups and product systems)?

\item
Is $E^\odot$ being faithful (that is, all $\rho_t$ are faithful) and strongly continuous in the new sense enough for showing existence of a strongly continuous right dilation (concluding the question of existence of faithful nondegenerate representations for product systems)?
\end{itemize}
We think that if the second and the third question have affirmative answers, then so does the first question in the case of strongly full and faithful product systems, and the proof will go exactly as we did for Theorem \ref{ecthm}. This would conclude the symmetric theory at least for strongly full and faithful strongly continuous (in the new sense) product systems of concrete von Neumann correspondences.

In Remark \ref{rem1st} we outlined already a possibility how to free this from the requirement on the product systems to be strongly full and faithful. As said in both Remarks, \ref{rem1st} and \ref{rem2nd}, these are projects for future work.

\section{Miscellaneous}

An immediate consequence of Theorem \ref{cpairedthm} and the results in \cite[Appendix B]{Ske16} is the following:

\bthm \label{cocthm}
If two strongly continuous \nbd{E_0}semigroups $\vt^1$ and $\vt^2$ on $\cB$ can be paired via strongly continuous pairings with the same strongly continuous \nbd{E_0}semigroup $\vt'$ on $\cB'$, then $\vt^1$ and $\vt^2$ are cocycle equivalent via a strongly continuous cocycle.
\ethm

\proof
$\vt^1$ and $\vt^2$ have isomorphic strongly continuous product systems, so the statement is granted by \cite[Theorem 12.3]{Ske16}.\qed

\brem
If $\vt^1$ and $\vt^2$, no matter whether strongly continuous or not, are paired via not necessarily strongly continuous pairings, then we only get a (not necessarily strongly continuous) isomorphism up to a multiplier. Imitating how in \cite[Theorem 6.5]{Ske16}, for isomorphic product systems, the cocycle $\eu_t\in\sB(G)$ is defined using the isomorphism $u_t$ (and the left dilations of $\vt^1$ and $\vt^2$), in the case here we get a family of unitaries $\eu_t$ (intertwining $\vt^1$ and $\vt^2$, as before) that fulfills the cocycle condition modulo the multiplier $m$. Like for $U_t$, the multiplier can be discussed away, only if it is trivial. (This happens, for instance, if $\eu_t$ and $\vt^1$, hence also $\vt^2$, are required strongly continuous.) We do not give more details, nor formulate a theorem.
\erem

\bulletline
Clearly, the proof of Theorem \ref{pairedthm} goes through without any change, if we restrict to discrete times $t\in\N_0$. Since the multipliers over $\N_0$ are trivial, we get the following result.

\bthm
Two discrete \nbd{E_0}semigroups $\vt$ and $\vt'$ on $\cB$ and $\cB'$, respectively, can be paired if and only if their product systems are commutants of each other.
\ethm

\proof
We only have to prove that a multiplier $m$ over $\N_0$ is trivial. Putting in the multiplier condition in \eqref{mult} $r=1$, $s=j$, and $t=k$, we get $m(1+j,k)=m(j,k)\frac{m(1,j+k)}{m(1,j)}$. For fixed $k$ and given $m(1,i)$, this recursion for $m(j,k)$ is resolved by
\beqn{
m(j,k)
~=~
\frac{m(1,j-1+k)}{m(1,j-1)}~\frac{m(1,j-2+k)}{m(1,j-2)}\,\ldots\,\frac{m(1,k)}{m(1,0)}m(0,k)
}\eeqn
($j\ge1$). Therefore, putting
\beqn{
f(i)
~:=~
\frac{m(0,0)}{m(1,0)\ldots m(1,i-1)}
}\eeqn
for $i\ge1$ and $f(0)=1$, taking also into account that $m(0,i)=m(0,0)=m(i,0)$ is constant, we get $m(j,k)=\frac{f(j)f(k)}{f(j+k)}$.\qed

\lf
We formulated this theorem in terms of semigroups (over $\N_0$) and product systems (over $\N_0$), because for this setting it is clear what we mean by the product system of a semigroup. But a product system over $\N_0$ is (isomorphic to) $\bfam{E_1^{\odot n}}_{n\in\N_0}$. By the correspondence of a unital endomorphism $\vt$ we mean $E_1$ from the product system of the the whole semigroup $\bfam{\vt^n}_{n\in\N_0}$. Now we can formulate a single mapping version of the theorem.

\bcor
For a unital normal endomorphism $\vt$ of $\cB$ and a unital normal endomorphism $\vt'$ of $\cB'$ there exists an automorphism $\alpha$ of $\sB(G)\supset\cB,\cB'$ such that $\alpha^{-1}\upharpoonright\cB=\vt$ and $\alpha\upharpoonright\cB'=\vt'$ if and only if the correspondences of $\vt$ and $\vt'$ are commutants of each other.
\ecor

Both the theorem and its corollary may equally well be proved directly (showing first the corollary directly and, then, extending the corollary to the whole semigroups generated by $\vt$ and $\vt'$ and the their respective correspondences).

{
\setlength{\baselineskip}{2.5ex}


\begin{thebibliography}{HKK04}

\bibitem[AHK78]{AlHK78}
S.~Albeverio and R.~Hoegh-Krohn, \emph{{Frobenius theory for positive maps of
  von Neumann algebras}}, Commun.\ Math.\ Phys. \textbf{64} (1978), 83--94.

\bibitem[AK92]{ArKi92}
W.~Arveson and A.~Kishimoto, \emph{{A note on extensions of semigroups of
  $*$--en\-do\-mor\-phisms}}, Proc.\ Amer.\ Math.\ Soc. \textbf{116} (1992),
  769--774.

\bibitem[Ale04]{Ale04}
A.~Alevras, \emph{{One parameter semigroups of endomorphisms of factors of type
  II$_1$}}, J.\ Operator Theory \textbf{51} (2004), 161--179.

\bibitem[Ara70]{Ara70}
H.~Araki, \emph{{Factorizable representations of current algebra}}, Publ.\
  Res.\ Inst.\ Math.\ Sci. \textbf{5} (1970), 361--422.

\bibitem[Arv69]{Arv69}
W.~Arveson, \emph{{Subalgebras of $C^*$--algebras}}, Acta Math. \textbf{123}
  (1969), 141--224.

\bibitem[Arv89]{Arv89}
\bysame, \emph{{Continuous analogues of Fock space}}, Mem.\ Amer.\ Math.\ Soc.,
  no. 409, American Mathematical Society, 1989.

\bibitem[Arv90]{Arv90}
\bysame, \emph{{Continuous analogues of Fock space IV: essential states}}, Acta
  Math. \textbf{164} (1990), 265--300.

\bibitem[Arv03]{Arv03}
\bysame, \emph{{Noncommutative dynamics and $E$--semigroups}}, Monographs in
  Mathematics, Springer, 2003.

\bibitem[Arv06]{Arv06}
\bysame, \emph{{On the existence of $E_0$--semigroups}}, Infin.\ Dimens.\
  Anal.\ Quantum Probab.\ Relat.\ Top. \textbf{9} (2006), 315--320.

\bibitem[BDH88]{BDH88}
M.~Baillet, Y.~Denizeau, and J.-F. Havet, \emph{{Indice d'une esperance
  conditionnelle}}, Compositio Math. \textbf{66} (1988), 199--236.

\bibitem[Bha96]{Bha96}
B.V.R. Bhat, \emph{{An index theory for quantum dynamical semigroups}}, Trans.\
  Amer.\ Math.\ Soc. \textbf{348} (1996), 561--583.

\bibitem[BISS14]{BISS14}
P.~Bikram, M.~Izumi, R.~Srinivasan, and V.~S. Sunder, \emph{{On extendability
  of endomorphisms and of {$E_0$}-semigroups on factors}}, Kyushu J. Math.
  \textbf{68} (2014), 165--179.

\bibitem[BM10]{BhMu10}
B.V.R. Bhat and M.~Mukherjee, \emph{{Inclusion systems and amalgamated products
  of product systems}}, Infin.\ Dimens.\ Anal.\ Quantum Probab.\ Relat.\ Top.
  \textbf{13} (2010), 1--26, (ar\-Xiv: 0907.0095v1).

\bibitem[BS00]{BhSk00}
B.V.R. Bhat and M.~Skeide, \emph{{Tensor product systems of Hilbert modules and
  dilations of completely positive semigroups}}, Infin.\ Dimens.\ Anal.\
  Quantum Probab.\ Relat.\ Top. \textbf{3} (2000), 519--575, (Rome,
  Volterra-Pre\-print 1999/0370).

\bibitem[Dav76]{Dav76}
E.B. Davies, \emph{{Quantum theory of open systems}}, Academic Press, 1976.

\bibitem[Dix77]{Dix77}
J.~Dixmier, \emph{{$C^\ast$--Algebras}}, North Holland Publishing Company,
  1977.

\bibitem[Fow02]{Fow02}
N.J. Fowler, \emph{{Discrete product systems of Hilbert bimodules}}, Pac.\ J.\
  Math. \textbf{204} (2002), 335--375.

\bibitem[Goh04]{Goh04}
R.~Gohm, \emph{{Noncommutative stationary processes}}, Lect.\ Notes Math., no.
  1839, Springer, 2004.

\bibitem[GS05]{GoSk05}
R.~Gohm and M.~Skeide, \emph{{Constructing extensions of CP-maps via tensor
  dilations with the help of von Neumann modules}}, Infin.\ Dimens.\ Anal.\
  Quantum Probab.\ Relat.\ Top. \textbf{8} (2005), 291--305, (ar\-Xiv:
  math.OA/0311110).

\bibitem[HKK04]{HKK04p}
J.~Hellmich, C.~K\"ostler, and B.~K\"ummerer, \emph{{Noncommutative continuous
  Bernoulli shifts}}, Pre\-print, ar\-Xiv: \newline math.OA/0411565, 2004.

\bibitem[Lie09]{Lie09}
V.~Liebscher, \emph{{Random sets and invariants for (type II) continuous tensor
  product systems of Hilbert spaces}}, Mem.\ Amer.\ Math.\ Soc., no. 930,
  American Mathematical Society, 2009, (ar\-Xiv: math.PR/0306365).

\bibitem[MS02]{MuSo02}
P.S. Muhly and B.~Solel, \emph{{Quantum Markov processes (correspondences and
  dilations)}}, Int.\ J.\ Math. \textbf{51} (2002), 863--906, (ar\-Xiv:
  math.OA/0203193).

\bibitem[MS05]{MuSo05}
\bysame, \emph{{~~~Duality of $W^*$-correspondences and applications}}, Quantum
  Probability and Infinite Dimensional Analysis --- From Foundations to
  Applications (M.~Sch\"urmann and U.~Franz, eds.), Quantum Probability and
  White Noise Analysis, no. XVIII, World Scientific, 2005, pp.~396--414.

\bibitem[MS13]{MaSr13}
O.T. Margetts and R.~Srinivasan, \emph{{Invariants for $E_0$--semigroups on
  II$_1$ factors}}, Commun.\ Math.\ Phys. \textbf{323} (2013), 1155--1184.

\bibitem[MS17]{MaSr17}
\bysame, \emph{{\hfill Non-cocycle-conjugate $E_0$--semigroups on factors}},
  Publ.\ Res.\ Inst.\ Math.\ Sci. \textbf{53} (2017), 299--336, (ar\-Xiv:
  1404.5934v2).

\bibitem[MSS06]{MSS06}
P.S. Muhly, M.~Skeide, and B.~Solel, \emph{{Representations of $\sB^a(E)$}},
  Infin.\ Dimens.\ Anal.\ Quantum Probab.\ Relat.\ Top. \textbf{9} (2006),
  47--66, (ar\-Xiv: math.OA/0410607).

\bibitem[Pas73]{Pas73}
W.L. Paschke, \emph{{Inner product modules over $B^*$--algebras}}, Trans.\
  Amer.\ Math.\ Soc. \textbf{182} (1973), 443--468.

\bibitem[Pow87]{Pow87}
R.T. Powers, \emph{{A non-spatial continuous semigroup of $*$--endomorphisms of
  $\sB(\eH)$}}, Publ.\ Res.\ Inst.\ Math.\ Sci. \textbf{23} (1987), 1053--1069.

\bibitem[Pow88]{Pow88}
\bysame, \emph{{An index theory for semigroups of $*$--endomorphisms of
  $\eB(\eH)$}}, Can.\ Jour.\ Math. \textbf{40} (1988), 86--114.

\bibitem[PR89]{PoRo89}
R.T. Powers and D.~Robinson, \emph{{An index for continuous semigroups of
  $*$--en\-do\-mor\-phisms of $\sB(H)$}}, J.\ Funct.\ Anal. \textbf{84} (1989),
  85--96.

\bibitem[PS72]{PaSchm72}
K.R. Parthasarathy and K.~Schmidt, \emph{{Positive definite kernels, continuous
  tensor products, and central limit theorems of probability theory}}, Lect.\
  Notes Math., no. 272, Springer, 1972.

\bibitem[Rie74]{Rie74a}
M.A. Rieffel, \emph{{Morita equivalence for $C^*$--algebras and
  $W^*$--algebras}}, J.\ Pure Appl.\ Algebra \textbf{5} (1974), 51--96.

\bibitem[Sch93]{MSchue93}
M.~Sch\"urmann, \emph{{White noise on bialgebras}}, Lect.\ Notes Math., no.
  1544, Springer, 1993.

\bibitem[Ske00]{Ske00b}
M.~Skeide, \emph{{Generalized matrix $C^*$--algebras and representations of
  Hilbert modules}}, Mathematical Proceedings of the Royal Irish Academy
  \textbf{100A} (2000), 11--38, (Cott\-bus, Rei\-he Mathe\-ma\-tik 1997/M-13).

\bibitem[Ske01]{Ske01}
\bysame, \emph{{Hilbert modules and applications in quantum probability}},
  Ha\-bi\-li\-ta\-tions\-schrift, Cottbus, 2001, available at:
  \\{\footnotesize\url{http://web.unimol.it/skeide/_MS/downloads/habil.pdf}}.

\bibitem[Ske02]{Ske02}
\bysame, \emph{{Dilations, product systems and weak dilations}}, Math.\ Notes
  \textbf{71} (2002), 836--843.

\bibitem[Ske03a]{Ske03c}
\bysame, \emph{{Commutants of von Neumann modules, representations of
  $\sB^a(E)$ and other topics related to product systems of Hilbert modules}},
  Advances in quantum dynamics (G.L. Price, B.M. Baker, P.E.T. Jorgensen, and
  P.S. Muhly, eds.), Contemporary Mathematics, no. 335, American Mathematical
  Society, 2003, (Preprint, Cottbus 2002, ar\-Xiv: math.OA/0308231),
  pp.~253--262.

\bibitem[Ske03b]{Ske03b}
\bysame, \emph{{Dilation theory and continuous tensor product systems of
  Hilbert modules}}, Quantum Probability and Infinite Dimensional Analysis
  (W.~Freudenberg, ed.), Quantum Probability and White Noise Analysis, no.~XV,
  World Scientific, 2003, (Preprint, Cottbus 2001), pp.~215--242.

\bibitem[Ske05]{Ske05a}
\bysame, \emph{{Three ways to representations of $\sB^a(E)$}}, Quantum
  Probability and Infinite Dimensional Analysis --- From Foundations to
  Applications (M.~Sch\"urmann and U.~Franz, eds.), Quantum Probability and
  White Noise Analysis, no. XVIII, World Scientific, 2005, (ar\-Xiv:
  math.OA/0404557), pp.~504--517.

\bibitem[Ske06a]{Ske06}
\bysame, \emph{{~A simple proof of the fundamental theorem about Arveson
  systems}}, Infin.\ Dimens.\ Anal.\ Quantum Probab.\ Relat.\ Top. \textbf{9}
  (2006), 305--314, (ar\-Xiv: math.OA/0602014).

\bibitem[Ske06b]{Ske06b}
\bysame, \emph{{Commutants of von Neumann correspondences and duality of
  Eilen\-berg-Watts theorems by Rieffel and by Blecher}}, Banach Center
  Publications \textbf{73} (2006), 391--408, (ar\-Xiv: math.OA/0502241).

\bibitem[Ske06c]{Ske06a}
\bysame, \emph{{Existence of $E_0$--semigroups for Arveson systems: Making two
  proofs into one}}, Infin.\ Dimens.\ Anal.\ Quantum Probab.\ Relat.\ Top.
  \textbf{9} (2006), 373--378, (ar\-Xiv: math.OA/0605480).

\bibitem[Ske06d]{Ske06d}
\bysame, \emph{{The index of (white) noises and their product systems}},
  Infin.\ Dimens.\ Anal.\ Quantum Probab.\ Relat.\ Top. \textbf{9} (2006),
  617--655, (Rome, Volterra-Pre\-print 2001/0458, ar\-Xiv: math.OA/0601228).

\bibitem[Ske08a]{Ske08}
\bysame, \emph{{Isometric dilations of representations of product systems via
  commutants}}, Int.\ J.\ Math. \textbf{19} (2008), 521--539, (ar\-Xiv:
  math.OA/0602459).

\bibitem[Ske08b]{Ske08a}
\bysame, \emph{{Product systems; a survey with commutants in view}}, Quantum
  Stochastics and Information (V.P. Belavkin and M.~Guta, eds.), World
  Scientific, 2008, (Preferable version: ar\-Xiv: 0709.0915v1!), pp.~47--86.

\bibitem[Ske09]{Ske09}
\bysame, \emph{{Unit vectors, Morita equivalence and endomorphisms}}, Publ.\
  Res.\ Inst.\ Math.\ Sci. \textbf{45} (2009), 475--518, (ar\-Xiv:
  math.OA/0412231v5 (Version 5)).

\bibitem[Ske11]{Ske11a}
\bysame, \emph{{Nondegenerate representations of continuous product systems}},
  J.\ Operator Theory \textbf{65} (2011), 71--85, (ar\-Xiv: math.OA/0607362).

\bibitem[Ske16]{Ske16}
\bysame, \emph{{Classification of $E_0$--semigroups by product systems}}, Mem.\
  Amer.\ Math.\ Soc., no. 1137, American Mathematical Society, 2016,
  electronically Oct 2015. (ar\-Xiv: 0901.1798v4).

\bibitem[Ske22]{Ske22b}
\bysame, \emph{{Hilbert von Neumann modules versus concrete von Neumann
  modules}}, Infinite Dimensional Analysis, Quantum Probability and
  Applications, ICQPRT 2021 (L.~Accardi, F.~Mukhamedov, and A.~Al~Rawashdeh,
  eds.), Springer Proceedings in Mathematics \& Statistics, no. 390, 2022,
  (ar\-Xiv: 1205.6413v2), pp.~169--182.

\bibitem[SS09]{ShaSo09}
O.M. Shalit and B.~Solel, \emph{{Subproduct systems}}, Documenta Math.
  \textbf{14} (2009), 801--868, (ar\-Xiv: 0901.1422v2).

\bibitem[SS23]{ShaSk23}
O.M. Shalit and M.~Skeide, \emph{{CP-Semigroups and dilations, subproduct
  systems and superproduct systems: The multi-parameter case and beyond}},
  Dissertationes Math. \textbf{585} (2023), 1--233, (ar\-Xiv: 2003.05166v3).

\bibitem[Wig39]{Wig39}
E.P. Wigner, \emph{{On unitary representations of the inhomogeneous Lorenz
  group}}, Ann.\ of Math. \textbf{40} (1939), 149--204.

\end{thebibliography}

\newcommand{\Swap}[2]{#2#1}\newcommand{\Sort}[1]{}
\providecommand{\bysame}{\leavevmode\hbox to3em{\hrulefill}\thinspace}
\providecommand{\MR}{\relax\ifhmode\unskip\space\fi MR }
\providecommand{\MRhref}[2]{%
  \href{http://www.ams.org/mathscinet-getitem?mr=#1}{#2}
}
\providecommand{\href}[2]{#2}

}

\lf\noindent
Michael Skeide: \it{Dipartimento di Economia, Universit\`{a} degli Studi del Molise, Via de Sanctis, 86100 Campobasso, Italy},
E-mail: \href{mailto:skeide@unimol.it}{\tt{skeide@unimol.it}},\\
Homepage: \url{http://web.unimol.it/skeide/}


\end{document}